\documentclass[12pt]{amsart}

\renewcommand{\Re}{\operatorname{Re}}

\newcommand{\sech}{\operatorname{sech}}
\newcommand{\defeq}{\stackrel{\rm{def}}{=}}

\usepackage{amssymb}
\usepackage{amsthm}
\usepackage{amsxtra}

\usepackage{fancyhdr}
\usepackage{nccmath}
\usepackage{graphicx}
\usepackage{color}
\usepackage{amsmath}
\usepackage{listings} 
\usepackage{float}
\usepackage{cases} 
\usepackage{threeparttable}
\usepackage{dcolumn}
\usepackage{multirow}
\usepackage{booktabs}

\usepackage{caption}
\usepackage{subcaption}

\newtheorem{theorem}{Theorem}[section]

\newtheorem{proposition}{Proposition}[section]

\theoremstyle{remark}
\newtheorem{remark}[proposition]{Remark}

\numberwithin{equation}{section}

\ifx\pdfoutput\undefined
  \DeclareGraphicsExtensions{.pstex, .eps}
\else
  \ifx\pdfoutput\relax
    \DeclareGraphicsExtensions{.pstex, .eps}
  \else
    \ifnum\pdfoutput>0
      \DeclareGraphicsExtensions{.pdf}
    \else
      \DeclareGraphicsExtensions{.pstex, .eps}
    \fi
  \fi
\fi

\title[Conservative gBO schemes]{High order conservative schemes for the generalized Benjamin-Ono equation in the unbounded domain}

\author{Kai Yang}
\date{\today}
\address{Florida International University}

\begin{document}
\maketitle
\begin{abstract}
This paper proposes a new class of mass or energy conservative numerical schemes for the generalized Benjamin-Ono (BO) equation on the whole real line with arbitrarily high-order accuracy in time. The spatial discretization is achieved by the pseudo-spectral method with the rational basis functions, which can be implemented by the Fast Fourier transform (FFT) with the computational cost $\mathcal{O}( N\log(N))$. By reformulating the spatial discretized system into the different equivalent forms, either the spatial semi-discretized mass or energy can be preserved exactly under the continuous time flow. Combined with the symplectic Runge-Kutta, with or without the scalar auxiliary variable reformulation, the fully discrete energy or mass conservative scheme can be constructed with arbitrarily high-order temporal accuracy, respectively.
Our numerical results show the conservation of the proposed schemes, and also the superior accuracy and stability to the non-conservative (Leap-frog) scheme.
\end{abstract}


\section{Introduction}
This paper considers the numerical methods for solving the generalized Benjamin-Ono (gBO) equation
\begin{align}\label{E:gBO}
\begin{cases}
u_t=-(-\mathcal Hu_x+\frac{1}{m}u^m)_x, \quad x \in \mathbb{R}, \,\, t>0, \,\, m \in \mathbb{Z}^+ ,\\
u(x,0)=u_0,
\end{cases} 
\end{align}
where the Hilbert transform $\mathcal{H}$ is defined by
\begin{align}\label{D:HT}
\mathcal Hf(x)=\frac{1}{\pi}\mbox{p.v.}\int_{-\infty}^{\infty} \frac{f(y)}{x-y} dy,
\end{align}
or equivalently, $\widehat{\mathcal Hf}(\xi)=-i\mbox{sgn}(\xi)\hat{f}(\xi)$ on the Fourier frequency side. When $m=2$, it is the well-known Benjamin-Ono (BO) equation 
\begin{equation}\label{E:BO}
u_t -\mathcal H u_{xx} + u_x u =0, 
\end{equation}
derived by Benjamin \cite{benjamin_1967} in 1967 and Ono \cite{Ono_1975} in 1975. This equation \eqref{E:BO} models the one-dimensional waves in deep water. The BO equation is closely related to the Korteweg-de Vries (KdV) equation, where the Hilbert transform term $\mathcal{H}u_{xx}$ is replaced by $u_{xxx}$. The KdV equation models the one-dimensional shallow water waves. 
Both equations, BO and KdV, are completely integrable. For example, the Lax pair can be constructed as described in e.g., \cite{NakaAk}, \cite{KM_98} \cite{AFA_83}, \cite{FA_83}
and \cite{AC1991}. Other nonlinearities for the equation \eqref{E:gBO} are also considered. They are relevant in various other models of water waves, e.g., see \cite{davis_acrivos_1967}, \cite{ABFS1989}, \cite{BK_79} and \cite{KB2000}. When $m=3$, the equation \eqref{E:gBO} is typically referred to as the modified Benjamin-Ono (mBO) equation. When $m \geq 3$, the equation \eqref{E:gBO} is typically referred to as the generalized Benjamin-Ono (gBO) equation.
In general, the gBO equation \eqref{E:gBO} conserves the following three quantities
\begin{align}
&I[u(t)] \defeq \int u(x,t) dx =I[u_0];   \label{E:momentum}\\
&M[u(t)] \defeq \int [u(x,t)]^2 dx =M[u_0];  \label{E:mass} \\
&E[u(t)] \defeq \int \left[ \frac{1}{2} \left((\mathcal H \partial_x)^{\frac{1}{2}} u(x,t)\right)^2-\frac{1}{m(m+1)} \left(u(x,t)\right)^{m+1} \right] dx = E[u_0]. \label{E:energy}
\end{align} 
The first one is called the $L^1$-type integral, and 
the last two are often called mass and energy (Hamiltonian), respectively.

Besides its physical applications, the gBO equation is also interesting to study from  the mathematical point of view. The well-posedness theory for the Cauchy problem has been discussed initially in \cite{S1979} and \cite{KPV1994}. Futher improvements on the well-posedness questions were done in \cite{KT2003}, \cite{KK2003}, \cite{T2004}, \cite{BP2008}, \cite{MR2004a}, \cite{MR2004b}, \cite{BP2006}, \cite{V2009}.
We also mention that when $m\geq 3$, it is typically referred to as the $L^2$-critical and $L^2$-supercritical cases from the scaling invariance, respectively. In those cases, there may exist blow-up solutions. This was numerically observed in \cite{BK04} and our recent paper \cite{RWY2021}. Besides the blow-up solutions, there are still many open questions, such as the soliton stability and the dispersion limit. These kind of questions have been studied both numerically and analytically. Compared with the (generalized) KdV equation (e.g., \cite{BSS1987}, \cite{GK2007}, \cite{GK2012}, etc.), the gBO equation is less well studied (e.g., \cite{MPST93}, \cite{PD00} and review \cite{S2018}). 
Therefore, a stable, efficient and accurate numerical algorithm would facilitate the future study. 

Numerical investigations on the BO equation have been started some time ago. Related articles can be found in \cite{BK04}, \cite{DHKR16}, \cite{TM98} for the domain truncation approach; \cite{Bo87}, \cite{Boyd1990}, \cite{Weideman95} for the computation of the Hilbert transform on $\mathbb{R}$; \cite{JW92}, \cite{BX12} for the pseudo-spectral method with the rational basis functions; and \cite{BX11} for a comparison between the domain truncation and the pseudo-spectral method on $\mathbb{R}$. 
Despite some years of investigations,
there are still far less studies about numerical methods for the gBO equations than the gKdV equations. To our best knowledge, there are no results concerning the conservative schemes for the gBO equation on the whole real line $\mathbb{R}$ so far. On the other hand, the conservative schemes are always preferable in simulating the  PDE's with conserved quantities, especially for studying the long time solution behavior, since it generally possess good accuracy and stability.
One possible reason is the 
numerical approximation of the Hilbert transform on $\mathbb{R}$, which is not as well studied as on a finite domain. However, if considering the conventional domain truncation spatial discretization strategy (e.g., the finite difference or Fourier spectral methods), the Hilbert transform usually leads to a slow decaying function, and consequently,
to a relatively large domain truncation error.

The purpose of this paper is to construct the conservative schemes for the gBO equation \eqref{E:gBO} on the whole real line $\mathbb{R}$, with arbitrarily high order accuracy in time. 
The spatial discretization is achieved by the rational basis functions with the pseudo-spectral approach from \cite{JW92}. We prove that by reformulating into the different forms, and applying the Hermitian or anti-Hermitian properties of the resulting spatial semi-discretized system, either the spatial semi-discretized mass or energy will be preserved. For the temporal discretization, the Crank-Nicholson method with the conventional reformulation of the nonlinear potential term (e.g., see \cite{DM2008} for the nonlinear Schr\"odinger (NLS) equation case) will lead to the conservation of the mass and energy in the discrete time flow. Furthermore, the high order conservative scheme can be constructed from the scalar auxiliary variable (SAV) approach (see \cite{SXY2018}, \cite{SXY2019} and \cite{CWJ2021}, \cite{Yang2021} for applications to dispersive PDEs). By using the symplectic Runge-Kutta (SRK) method, the three invariant quantities \eqref{E:momentum}-\eqref{E:energy} (with proper modifications for energy \eqref{E:energy}) will be preserved exactly in the discrete time flow. However, due to the limitation of the spatial discretization, we can only conserve either the discrete mass or the discrete energy in the space-time fully discrete sense. 
In fact, this strategy is universal. By a similar space-time discretization, it is easy to construct the conservative schemes for the gKdV equations and the structure-preserving schemes (the discrete mass and energy are preserved exactly at the same time) for the NLS equations, as well as their high dimensional generalization by applying the tensor product. 
This will be useful in studying the long time behavior of the solutions for those equations, as well as the slow decaying solutions, since the traditional domain truncation strategy (e.g., \cite{LY2006}, \cite{YHL2013}, \cite{LY2016}, and \cite{GWWC2017}) requires large computational domain, and consequently, it results in large number of nodes in spatial discretization.

This paper is organized as follows. In Section 2, we introduce the pseudo-spectral spatial discretization strategy from the rational basis functions. Then, we define the discrete inner product with respect to the collocation points from such rational basis functions.
Finally, we give the mass-conservative or the energy-conservative spatial semi-discretized form of the gBO equation \eqref{E:gBO}. In Section 3, we first introduce the Crank-Nicholson types of temporal discretization. We show that the Crank-Nicholson method with its conventional modification on the nonlinear term, will preserve the mass and energy exactly in the discrete time flow.
Combining with the previous results in Section 2, we give two fully discretized schemes for the gBO equation, which conserve either the discrete mass or the discrete energy. Next, we consider the high-order conservative schemes achieved by the symplectic Runge-Kutta method with the SAV reformulation. From the classical argument (e.g. \cite{Cooper1987}, \cite{Sanz1988} and \cite{Yang2021}), we show that the reformulated system preserves the quantities \eqref{E:momentum}-\eqref{E:energy} exactly in the discrete time flow. Again, combined with the spatial discretization results in Section 2, we give two fully discretized schemes with high order temporal accuracy---one conserves the discrete mass and the other conserves the discrete energy. 
In Section 4, we illustrate the numerical examples. As a comparison, we also show the numerical results obtained from the non-conservative semi-implicit Leap-Frog scheme. Our numerical results show that the proposed schemes preserve the designate quantities based on the type of conservative scheme we choose. Compared with the non-conservative scheme, these conservative schemes also possess better accuracy. Additionally, the error from the temporal discretization decreases on the order as expected (second order for the IRK2 and Leap-Frog schemes, and fourth order for the IRK4 schemes). These results show the validity and efficiency of the numerical methods proposed.

\medskip
{\flushleft \bf Acknowledgment:} The author is partially supported by the NSF grant DMS-1927258 (PI: Svetlana Roudenko). The author is thankful for Dr. Roudenko's helpful discussion, reading and remarks on the paper.

\section{Spatial discretization}\label{S:SD}
In this section, we describe the rational basis functions in $\mathbb{R}$ used for the spatial discretization. The review of basis functions can be found in \cite{Christov82} and \cite{Weideman95}. One advantage of this discretization is that it can easily represent the Hilbert transform.
Then, we define the discrete inner product corresponding to the collocation points from the rational basis functions.
Finally, we introduce two types of the spatial discretization for the gBO equation \eqref{E:gBO}: one is mass-conservative and the other one is energy-conservative.

\subsection{Rational basis functions}
Consider the rational basis functions on the whole real line $\mathbb{R}$,
\begin{equation}\label{rational}
 u(x,t)=\sum_{k=-\infty}^{\infty} \hat{u}_{k}(t)\rho_k(x), \quad \rho_k(x)=\frac{(\alpha+ix)^k}{(\alpha-ix)^{k+1}},
\end{equation}
where $\alpha$ is a mapping parameter that we will describe later. In \cite{Christov82}, it is shown that $\{ \rho_k(x) \}_{k=-\infty}^{\infty}$ 
form a complete orthogonal basis in $L^2(-\infty, \infty)$ with the following orthogonality
\begin{align}\label{E:rho orthogonal}
\int_{-\infty}^\infty \rho_j(x) \overline{\rho_k(x)}dx = \begin{cases}
      \pi/\alpha, &  j=k\\
      0, & j\neq k.
    \end{cases}   \defeq \frac{\pi}{\alpha} \delta_{j,k} .
\end{align}
Therefore, we have $$\hat{u}_k(t)=\frac{\alpha}{\pi}\int_{-\infty}^{\infty} u(x,t)\rho_k(x)dx.$$
From the rational expansion \eqref{rational}, the Hilbert transform can be easily calculated \cite{Weideman95} by
\begin{align}\label{E:HT rational}
\mathcal{H}(u(t,x))=\sum_{k=-\infty}^{\infty} -i \hat{u}_k(t) \text{sgn}(k)\rho_k(x),
\end{align}
with $\mathrm{sgn}(k)=1$ when $k=0$.
Meanwhile, the derivatives of $u(x,t)$ can be computed by the relation
\begin{align}
u_x(x,t)=\sum_{k=-\infty}^{\infty} &\frac{i}{2\alpha}[k \hat{u}_{k-1}+(2k+1) \hat{u}_k+(k+1) \hat{u}_{k+1}]\rho_k(x), \label{ux}
\end{align}
and higher order derivatives could be done iteratively.

In numerical computations, a truncation of $N$-term interpolation function $I_Nu$ are used to approximate the function $u(x)$, i.e.,
$$ u(x,t)\approx I_Nu \defeq \hat{\mathbf{u}}^T \vec{\rho}:=\sum_{k=-N/2}^{N/2-1} \hat{u}_k(t)\rho_k(x),$$
where $\hat{\mathbf{u}}=(\hat{u}_{-N/2},\hat{u}_{-N/2+1}, \cdots, \hat{u}_{N/2-1})^T$ is the vector of the truncated coefficients, and $\vec{\rho}=(\rho_{-N/2}(x), \cdots, \rho_{N/2-1}(x))^T$ is the vector function of $\rho_k(x)$.
This leads to the sparse matrix forms
\begin{align}\label{E:rational basis derivative}
 u_x \approx [ \mathbf{S_1} \hat{\mathbf{u}}]^T \vec{\rho}, \quad u_{xx} \approx [\mathbf{S_2} \hat{\mathbf{u}}]^T \vec{\rho}, \quad \mathcal{H}u \approx [\mathbf{H}\hat{\mathbf{u}}]^T \vec{\rho},
\end{align}
where $\bf{S_{1}}$ is given in \eqref{ux} via the coefficients of $\lbrace \hat{u}_k \rbrace$, and $\mathbf{S_2}=\mathbf{S_1}\times \mathbf{S_1}$ by computing the derivatives iteratively from \eqref{ux}, and $\mathbf{H}=-i \mathrm{diag}(\mathrm{sgn}(-N/2+0.5), \cdots, \mathrm{sgn}(N/2-0.5))$ is the diagonal matrix representing the approximation of the Hilbert transform in \eqref{E:HT rational}.


From \eqref{ux} and \eqref{E:rational basis derivative}, it is easy to see that the matrices $\mathbf{S_1}$ and $\mathbf{H}$ are anti-Hermitian, and the matrix $\mathbf{S_2}$ is real and symmetric.

Now, consider the change of variable 
$$x=\alpha\tan \frac{\theta}{2}, \quad \mathrm{or} \,\, \mathrm{equivalently}, \quad e^{i\theta}=\frac{\alpha+ix}{\alpha-ix}, \quad -\pi\leq \theta \leq \pi,$$
and a spatial discretization $x_j=\alpha\tan \frac{\theta_j}{2}, \theta_j=jh, h=2\pi/N, j=-N/2, \cdots, N/2-1$, where the $\alpha$ is the mapping parameter indicating that $N/2$ collocation points are located in the interval $[-\alpha, \alpha]$. Notice that 
\begin{equation}\label{uj}
u(x_j)=\sum_{k=-N/2}^{N/2-1}\hat{u}_k \rho_k(x_j) \,\, \Rightarrow \,\,   u(x_j)(\alpha-ix_j)=\sum_{k=-N/2}^{N/2-1}\hat{u}_k e^{ik\theta_j},
\end{equation}
hence, the Fast Fourier transform (FFT) can be applied to obtain the coefficients $\hat{u}_k$.
We note that the above discretization in space is not uniform in $x$, but uniform in $\theta$, and the singularity at $x_{-N/2}=-\infty$ can be removed by imposing the boundary condition $u(-\infty)=0$, i.e., $u_{-N/2}=0$. 

We denote the matrix $\mathbf{F}$ to be the standard Fast Fourier transform (FFT) matrix with $\{ k\theta_j \}$, i.e.,
$$\mathbf{F}_{kj}= \frac{1}{N} e^{-ik  \theta_j}, \quad \mathbf{F}^{-1}_{jk}= e^{ik  \theta_j}, \quad -N/2<j,k<N/2-1.$$
Note that instead of writing explicitly, the matrices $\mathbf{F}$ and $\mathbf{F}^{-1}$ can be computed by FFT, (see, e.g., \cite[Chapter 2]{STL2011} and \cite[Chapter 3]{Tr2000}). Denote the diagonal matrix $\mathbf{P}=\text{diag}(\alpha-ix_{-N/2}, \cdots, \alpha-ix_{N/2-1})$ to be the weight matrix, which comes from \eqref{uj}. The coefficients $\hat{u}_k$ in the vector form can be represented by
\begin{align}\label{E:uhat}
\mathbf{\hat{u}}=\mathbf{FPu},
\end{align}
where $\mathbf{u}=(u_{-N/2},\cdots, u_{N/2-1})^T$ and $u_j=u(x_j)$.

Now, we define the discrete inner product with respect to the rational basis function. 
Denote the inner product between the two functions $u(x)$ and $v(x)$ on $\mathbb{R}$ by
$$\langle u(x),v(x) \rangle \defeq \int_{\mathbb{R}} u(x) \bar{v}(x) dx.$$
Recall that the interpolation function $I_N u$ is the approximation of $u(x) \approx I_Nu = \sum_{k=-N/2}^{N/2-1} \hat{u}_{k}\rho_k(x)$, then, the approximation of the inner product for functions $u$ and $v$ will be
\begin{align}\label{D:L2}
\langle u, v \rangle & \approx \int_{\mathbb{R}} I_Nu \overline{I_N v} \, dx = \frac{\pi}{\alpha}  \sum_{k=-N/2}^{N/2-1} \hat{u}_k \overline{\hat{v}}_k \\
& =\frac{\pi}{\alpha}\mathbf{(\overline{FPv})^T(FPu)}=\frac{\pi}{\alpha}\mathbf{\bar{v}^T \bar{P} \bar{F}^T F Pu}=\frac{\pi}{\alpha N}\mathbf{\bar{v}^T \bar{P}  Pu}, \nonumber
\end{align}
from the orthogonal property \eqref{E:rho orthogonal}, the relation \eqref{E:uhat}, and $\mathbf{\bar{F}^T}=\frac{1}{N} \mathbf{F^{-1}}$ (e.g., see \cite[Chapter 2]{STL2011}). Denote the diagonal matrix $\mathbf{W}=\mathbf{P\bar{P}}=\mathrm{diag}(\alpha^2+x^2_{-N/2}, \cdots, \alpha^2+x^2_{N/2-1})$ to be the product of the two diagonal matrices $\mathbf{P}$ and $\mathbf{\bar{P}}$. According to \eqref{D:L2}, we can define the discrete inner product with respect to the collocation points $\lbrace x_j \rbrace$ from the rational basis function as follows:
\begin{align}\label{D:L2 discrete}
\langle \mathbf{u}, \mathbf{v} \rangle_h \defeq \frac{\pi}{N \alpha} \mathbf{\bar{v}^TWu}=\frac{\pi}{N \alpha} \sum_{j=-N/2}^{N/2-1} w_j u_j \bar{v}_j,
\end{align}  
where $w_j=\alpha^2+x^2_{j}$ can be considered as the weights for the quadrature.

\subsection{Conservative spatial discretization}
To discuss the conservative spatial discretizations, we first define the spatial semi-discretized $L^1$-type integral, mass and energy from \eqref{E:momentum}-\eqref{E:energy}. Let $\mathbf{1}=(1,1, \cdots ,1)^T$ be the $N \times 1$ vector. For simplicity, we also denote by $\mathbf{u^m}=(u_{-N/2}^m, \cdots u_{N/2-1}^m)^T$ to be the pointwise power of the vector $\mathbf{u}$. 
Then, the spatial semi-discretized $L^1$-type integral, mass and energy  are defined as follows
\begin{align}
&I_h=\langle \mathbf{u}, \mathbf{1} \rangle_h; \label{D:momentum s} \\
&M_h=\langle \mathbf{u}, \mathbf{u} \rangle_h; \label{D:mass s} \\
&E_h=\frac{1}{2} \langle \mathbf{P^{-1}F^{-1} HS_1 FP u}, \mathbf{u} \rangle_h -\frac{1}{m(m+1)} \langle \mathbf{u^m}, \mathbf{u} \rangle_h. \label{D:energy s}
\end{align}

It is easy to see that if $\mathbf{u} \in \mathbb{R}^N$, then $\frac{d}{dt}I_h=\langle \mathbf{u}_t, \mathbf{1} \rangle_h$, and $\frac{d}{dt}M_h=2\langle \mathbf{u}_t, \mathbf{u} \rangle_h$ from \eqref{D:L2 discrete}. 
We also note that 
\begin{align}\label{E:Denergy dt}
\frac{d}{dt}E_h& = \frac{1}{2} \left( \langle \mathbf{P^{-1}F^{-1} HS_1 FP u}_t, \mathbf{u} \rangle_h+\langle \mathbf{P^{-1}F^{-1} HS_1 FP u}, \mathbf{u}_t \rangle_h \right)  \\
& - \frac{1}{m(m+1)} \left(\langle (\mathbf{u^m})_t, \mathbf{u} \rangle_h +\langle \mathbf{u^m}, \mathbf{u}_t \rangle_h \right) \nonumber \\
&=\Re \left( \langle \mathbf{P^{-1}F^{-1} HS_1 FP u}, \mathbf{u}_t \rangle_h \right) -\frac{1}{m} \langle \mathbf{u^m}, \mathbf{u}_t \rangle_h,  \nonumber
\end{align}
since $\langle \mathbf{P^{-1}F^{-1} HS_1 FP u}, \mathbf{u}_t \rangle_h =\overline{\langle \mathbf{P^{-1}F^{-1} HS_1 FP u_t}, \mathbf{u} \rangle_h}$, and 
$$\langle \mathbf{(u^m)}_t, \mathbf{u} \rangle_h =m \langle \mathbf{u^{m-1}u}_t, \mathbf{u} \rangle_h =m\langle \mathbf{u^m}, \mathbf{u}_t \rangle_h$$ 
from \eqref{D:L2 discrete}.

Now, we have the following proposition for the spatial conservative discretizations.
\begin{proposition}\label{P:spatial c}
The following spatial semi-discretized equation to the gBO equation \eqref{E:gBO}
\begin{align}\label{E:gBO SEC}
\mathbf{u}_t=-\mathbf{P^{-1}F^{-1}S_1FP}(\mathbf{-P^{-1}F^{-1}HS_1PF u}+ \frac{1}{m}\mathbf{u^m})
\end{align}
conserves the spatial semi-discretized energy, i.e.,
\begin{align}\label{E:Eh dt}
\frac{d}{dt}E_h=0.
\end{align}
On the other hand, the following spatial semi-discretized equation to the gBO equation \eqref{E:gBO}
\begin{small}
\begin{align}\label{E:gBO SMC}
\mathbf{u}_t=\mathbf{P^{-1}F^{-1}HS_2FP u}- \frac{1}{m+1} \left( \mathbf{\mathrm{diag}(u^{m-1}) P^{-1}F^{-1}S_1FP u} +\mathbf{P^{-1}F^{-1}S_1FPu^m} \right)
\end{align}
\end{small}
conserves the spatial semi-discretized mass, i.e.,
\begin{align}\label{E:Mh dt}
\frac{d}{dt}M_h=0.
\end{align}
\end{proposition}

\begin{proof}
Putting the equation \eqref{E:gBO SEC} in \eqref{E:Denergy dt} yields
\begin{small}
\begin{align}\label{Eh identity}
\frac{d}{dt}E_h&=\Re \left( \langle \mathbf{P^{-1}F^{-1}HS_1FP u}- \frac{1}{m}\mathbf{u^m}, \mathbf{P^{-1}F^{-1}S_1FP}(\mathbf{P^{-1}F^{-1}HS_1FP u}- \frac{1}{m}\mathbf{u^m})  \rangle_h   \right) \\
&=\frac{\pi}{\alpha N} \Re \left[ \overline{ \left( \mathbf{P^{-1}F^{-1}S_1FP}(\mathbf{P^{-1}F^{-1}HS_1FP u}- \frac{1}{m}\mathbf{u^m}) \right) } ^{\mathbf{T}} \mathbf{W} \left( \mathbf{P^{-1}F^{-1}HS_1FP u}- \frac{1}{m}\mathbf{u^m} \right)  \right] \nonumber \\
&= \frac{\pi}{\alpha N}  \Re \left[ \overline{ \left(\mathbf{P^{-1}F^{-1}HS_1FP u}- \frac{1}{m}\mathbf{u^m}\right) }^{\mathbf{T}}  \left( \mathbf{\bar{P}F^{-1}S_1FP} \right)^{\mathbf{T}} \left( \mathbf{P^{-1}F^{-1}HS_1FP u}- \frac{1}{m}\mathbf{u^m} \right) \right]=0, \nonumber
\end{align}
\end{small}
since the matrix $\frac{1}{N} \mathbf{\bar{P}F^{-1}S_1FP}$ is anti-Hermitian.

Similarly, we first note that $\mathbf{W\mathrm{diag}(u^{m-1})}=\mathbf{\mathrm{diag}(u^{m-1})W}$, since both of them are diagonal matrices. Then,
\begin{align}\label{Mh identity}
\frac{d}{dt}M_h&=2\langle \mathbf{u_t}, \mathbf{u} \rangle =2\langle \mathbf{P^{-1}F^{-1}HS_2FP u}, \mathbf{u} \rangle_h  \\
-&\frac{2}{m+1} \langle \left( \mathbf{\mathrm{diag}(u^{m-1}) P^{-1}F^{-1}S_1FP u} +\mathbf{ P^{-1}F^{-1}S_1FPu^m } \right) , \mathbf{u} \rangle_h  \nonumber \\
&=0-\frac{2\pi}{(m+1) \alpha N}\left( (\mathbf{u^{m})^T \bar{P}F^{-1}S_1FP u}+ \mathbf{u^T \bar{P}F^{-1}S_1FP u^m}  \right)=0,\nonumber
\end{align}
from the similar decomposition to \eqref{Eh identity}, and the fact that the matrices $\frac{1}{N} \mathbf{\bar{P}F^{-1}HS_2FP}$ and $\frac{1}{N} \mathbf{\bar{P}F^{-1}S_1FP}$  are also anti-Hermitian. 
\end{proof}

\section{Temporal and full discretization}\label{S:TD}
In this section, we first discuss the temporal discretization, and then, the space-time full discretization of the gBO equation \eqref{E:gBO}. We start with the most commonly used Crank-Nicholson-type scheme. After that, we consider the high order conservative schemes. This is achieved by the symplectic Runge-Kutta method, such as the Gauss-Legendre Runge-Kutta method. When considering the energy conservation, the scalar auxillary variabel (SAV) approach from \cite{SXY2018} and \cite{CWJ2021} will be incorporated.
With this kind of approach, one can easily construct the conservative numerical scheme with arbitrarily high order accuracy in time. 
\subsection{Crank-Nicholson-type scheme}
We first introduce the notations. 
Assume that our simulation is on the finite time interval $t\in [0,T]$. Define $\tau$ to be the time step and $t_n=n\tau$ to be the time at the $n$th time step. Denote $u^n \approx u(x,t_n)$ to be the semi-discretization in time.
Denote $I^n=I[u^n]$, $M^n=M[u^n]$ and $E^n=E[u^n]$ to be the momentum, mass and energy from \eqref{E:momentum}-\eqref{E:energy} at time $t=t_n$.
For convenience, the {half-time step} is denoted as $u^{n+\frac{1}{2}}=\frac{1}{2}(u^n+u^{n+1})$ from the linear interpolation.
We also denote the full discretization by $u^n_j \approx u(x_j,t_n)$, and the column vector $\mathbf{u}^n \approx u(\mathbf{x},t_n)$. Now, we define the discrete $L^1$-type integral, mass and energy as follows:
\begin{align}
&I_h^n=\langle \mathbf{u}^n, \mathbf{1} \rangle_h; \label{D:momentum st} \\
&M_h^n=\langle \mathbf{u}^n, \mathbf{u}^n \rangle_h; \label{D:mass st} \\
&E_h^n=\frac{1}{2} \langle \mathbf{P^{-1}F^{-1} HS_1 FP} \mathbf{u}^n, \mathbf{u}^n \rangle_h -\frac{1}{m(m+1)} \langle (\mathbf{u}^n)^{\mathbf{m}}, \mathbf{u}^n \rangle_h. \label{D:energy st}
\end{align}
We propose the following Crank-Nicholson-type  {\it{mass-conservative}} scheme.
\begin{theorem}\label{T:CN MC}
The scheme
\begin{small}
\begin{align}\label{CN MC}
\frac{\mathbf{u}^{n+1}-\mathbf{u}^n}{\tau}&=\mathbf{P^{-1}F^{-1}HS_2FP u}^{n+\frac{1}{2}} \\
&- \frac{1}{m+1} \left( \mathrm{diag}\left((\mathbf{u}^{{n+\frac{1}{2}}})^{\mathbf{m-1}} \right) \mathbf{P^{-1}F^{-1}S_1FP} \mathbf{u}^{n+\frac{1}{2}} +\mathbf{P^{-1}F^{-1}S_1FP}(\mathbf{u}^{{n+\frac{1}{2}}})^{\mathbf{m}} \right) \nonumber
\end{align}
\end{small}
conserves the discrete mass \eqref{D:mass st} exactly in time, i.e.,
$$M_h^{n+1}=M_h^{n}.$$
\end{theorem}
\begin{proof}
The proof is straightforward. Equipping the equation \eqref{CN MC} with the discrete inner product \eqref{D:L2 discrete} with the vector $\mathbf{u}^{n+\frac{1}{2}}$, and using the identity \eqref{Mh identity} in Proposition \ref{P:spatial c} yields the result.
\end{proof}

For the Crank-Nicholson-type {\it{energy-conservative}} scheme, by modifying of the nonlinear term, we have the following theorem.
\begin{theorem}\label{T:CN EC}
The scheme
\begin{align}\label{CN EC}
\frac{\mathbf{u}^{n+1}-\mathbf{u}^n}{\tau}&=\mathbf{P^{-1}F^{-1}S_1FP} \Big[  \mathbf{P^{-1}F^{-1}HS_1FP u}^{n+\frac{1}{2}} \\
&- \frac{1}{m(m+1)}  \mathrm{diag}\left( \frac{(\mathbf{u}^{n+1})^{\mathbf{m+1}}-(\mathbf{u}^{n})^{\mathbf{m+1}}}{(\mathbf{u}^{n+1})^{\mathbf{2}}-(\mathbf{u}^{n})^{\mathbf{2}} } \right) \mathbf{u}^{n+\frac{1}{2}} \Big] \nonumber
\end{align}
conserves the discrete energy \eqref{D:energy st} exactly in time, i.e.,
$$E_h^{n+1}=E_h^{n}.$$
\end{theorem}
\begin{proof}
Following the same idea as in Theorem \ref{T:CN MC}, we equip the equation \eqref{CN EC} with the discrete inner product \eqref{D:L2 discrete} with the vector 
$$\mathbf{P^{-1}F^{-1}HS_1FPu}^{n+\frac{1}{2}}- \frac{1}{m(m+1)}  \mathrm{diag}\left( \frac{(\mathbf{u}^{n+1})^{\mathbf{m+1}}-(\mathbf{u}^{n})^{\mathbf{m+1}}}{(\mathbf{u}^{n+1})^{\mathbf{2}}-(\mathbf{u}^{n})^{\mathbf{2}} } \right) \mathbf{u}^{n+\frac{1}{2}}. $$ 
Then, using the identity \eqref{Eh identity} in Proposition \ref{P:spatial c} yields the result.
\end{proof}

The construction for the conservative schemes in Theorem \ref{T:CN MC} and Theorem \ref{T:CN EC} are standard. If we only consider the semi-discretization in time, the scheme \eqref{CN MC} is the midpoint rule, or the implict 2nd order Runge-Kutta method (IRK2), which is also known as the symplectic (quadratic preserving) Runge-Kutta method. We split the potential into the form 
$$\frac{1}{m}(u^m)_x= \frac{1}{m+1}[u^{m-1}u_x +(u^m)_x]$$
 for the purpose of creating the symmetry for the mass conservations in the spatial discretization, which we discussed in the previous section. For the energy conservation, we need to reformulate the potential part, which is widely used in  literature, see e.g., \cite{DM2002}, \cite{MRRY2021} and \cite{LRY2021} for the NLS case. It is easy to see that if we only consider the temporal semi-discretization (assuming the spatial variable is continuous), then the scheme \eqref{CN EC} will conserve both the mass and energy in the discrete time flow.

We next discuss the numerical schemes with higher order temporal accuracy. For simplicity and conciseness, we only consider the semi-discretization in time. The space-time full discretization results can be easily generalized together with the results from Section \ref{S:SD}.
\subsection{High order conservative schemes}
The high order temporal conservative schemes can be achieved by the symplectic Runge-Kutta (SRK) method. 
We first briefly review the RK method before showing our results. Consider the problem
\begin{align}\label{E:RK}
u_t=f(u).
\end{align}
From the time $t=t_n$ to $t=t_{n+1}$, let $b_i$, $a_{ij}(i,j=1,\cdots s)$ be real numbers, and $c_i=\sum_{j=1}^s a_{ij}$ be the collocation points. Denote the intermediate values $U_i$ to be the solution satisfying \eqref{E:RK} at the intermediate time $t^i=t_n+\tau c_i$. Then, the intermediate values $U_i$'s are calculated by
\begin{align}\label{Ui}
U_i=u^n+\tau \sum_{j=1}^s a_{ij}f_j,
\end{align}
where $f_i=f(U_i)$. The solution $u^{n+1}$ is updated by 
\begin{align}\label{u n+1}
u^{n+1}=u^n+\tau \sum_{j=i}^s b_if_i.
\end{align}

We usually write the coefficients $\mathbf{A}=(a_{ij})$, $\mathbf{b}=(b_1, b_2,\cdots,b_s)$ and $\mathbf{c}=(c_1, c_2,\cdots,c_s)^T$ in the Butcher's Tableaus (\cite{Butcher1964}):
\[
\begin{array}
{c|c}
\mathbf{c}&
\mathbf{A}\\
\hline
& \mathbf{b}
\end{array}.
\]

For example, we list two commonly used Runge-Kutta methods in the Butcher's Tableaus in Table \ref{T:RK Butcher}. They are the $s$-stage Runge-Kutta methods with $s=1,2$, respectively. These methods are coming from the Gaussian-Legendre quadrature, known as the IRK2 and IRK4 methods, since the temporal accuracy is on the order of $2$ and $4$, respectively. We use these methods in our numerical simulations in the next section. There are many other types of Runge-Kutta methods as well, we refer the interested reader to \cite{BBM1972}, \cite{Cooper1987}, \cite{Sanz1988}, \cite{Geng1993} and \cite{SA1991}.
\begin{table}[h]
\renewcommand\arraystretch{1.2}
\centering
\begin{subtable}[t]{0.1\textwidth}
\[
\begin{array}{c|c}

\frac{1}{2}   &   \frac{1}{2}   \\ \hline
    &   1   \\ 
\end{array}
\]
\caption{IRK2}
\end{subtable}
    \hfil
\begin{subtable}[t]{0.3\textwidth}
\[
\begin{array}{c|cc}
   
\frac{1}{2}-\frac{1}{6}\sqrt{3}   &   \frac{1}{4}  & \frac{1}{4}-\frac{1}{6}\sqrt{3}  \\ 
\frac{1}{2}+\frac{1}{6}\sqrt{3}   &   \frac{1}{4}+\frac{1}{6}\sqrt{3}  & \frac{1}{4}  \\ \hline
    &   \frac{1}{2} & \frac{1}{2}   \\ 
\end{array}
\]
\caption{IRK4}
\end{subtable}
    \hfil
\caption{\label{T:RK Butcher}Butcher's Tableaus for the $s$-stage Gaussian-Legendre collocation Runge-Kutta methods with $s=1,2$.}
\end{table}

We prove the following theorem for the mass-conservative scheme.
\begin{theorem}\label{T:SRK MC}
The $s$-stage symplectic (quadratic preserving) Runge-Kutta method, which satisfies 
\begin{align}\label{E:SRK ab}
b_ia_{ij}+b_ja_{ji}=b_ib_j, \qquad \mbox{for} \quad i,j=1,\cdots,s,
\end{align} 
conserves the discrete mass exactly in time for the spatial discretized gBO equation \eqref{E:gBO SMC}, i.e.,
$$M_h^{n+1}=M_h^{n}.$$
\end{theorem}
\begin{proof}
The proof is standard. From the standard RK theory (e.g., \cite{Cooper1987}, \cite{Sanz1988}), we can show that the temporal semi-discretized scheme conserves the mass exactly in the discrete time flow. 
Indeed, the RK theory shows that 
$$M^{n+1}-M^n= 2\tau \sum_{i=1}^s b_i \langle U_i, f(U_i)  \rangle+ \tau^2 \sum_{i,j=1}^s (b_ia_{ij}+b_ja_{ji}-b_ib_j) \langle f(U_i), f(U_j)  \rangle=0,  $$
since $\langle U_i, f(U_i)  \rangle=0$ by putting $f(U)$ in the form of \eqref{E:gBO}. 

When considering the space-time full discratization, we have 
$$M^{n+1}_h-M^n_h= 2\tau \sum_{i=1}^s b_i \langle \mathbf{U}_i, f(\mathbf{U}_i)  \rangle_h+ \tau^2 \sum_{i,j=1}^s (b_ia_{ij}+b_ja_{ji}-b_ib_j) \langle f(\mathbf{U}_i), f(\mathbf{U}_j)  \rangle_h=0,  $$
where the vector $\mathbf{U}_i$ is the discretized version of the intermediate value $U_i$ for $i=1, \cdots, s$, and  $\langle \mathbf{U}_i, f(\mathbf{U}_i)  \rangle_h=0 $ by using the relation \eqref{Mh identity}.

\end{proof}

The symplectic Runge-Kutta method cannot preserve the discrete energy. In order to construct the energy-preserving scheme, we need to reformulate the potential term in the same idea as in \eqref{CN EC}. This is achieved by using the scalar auxillary approach from \cite{SXY2018} and \cite{CWJ2021}. We reformulate the equation \eqref{E:gBO} into an equivalent system as follows:
\begin{align}\label{E:gBO sav}
\begin{cases}
u_t=-\left( -\mathcal H u_{x}+\dfrac{1}{m} \frac{u^m v}{\sqrt{(u^m,u)+C_0}} \right)_x,\\
v_t=\frac{m+1}{2\sqrt{\langle u^m,u \rangle+C_0}} \langle u^m ,u_t \rangle,
\end{cases}
\end{align}
with the initial condition 
\begin{align*}
u(x,0)=u_0, \qquad v_0=\sqrt{\langle u_0^m,u_0 \rangle+C_0}.
\end{align*}
Then, the energy to the system \eqref{E:gBO} is modified into the equivalent form
\begin{align}\label{E:energy sav}
E[u(t), v(t)] \defeq \frac{1}{2} \langle \mathcal{H}u_{x},u \rangle -\frac{1}{m(m+1)} (v^2 -C_0) \equiv E[u_0, v_0].
\end{align}
Here, we slightly nabuse the notation $E[u,v]$ to represent the modified energy for convenience, since it is equivalent to the energy $E[u]$ in \eqref{E:energy} in the continuous sense.
The $C_0$ is a constant to make sure that the term $\langle u^m, u\rangle +C_0$ is  positive for all time $t \in [0,T]$. In the actual computation, the $C_0$ is adjustable during the time evolution, and thus, we only need to choose the constant $C_0$ such that the term $\langle u^m, u\rangle +C_0>0$ in the time interval $t\in [t_n,t_{n+1}]$. This is easily fulfilled, since we only consider the solution smooth in time. We will discuss the $C_0$ adjustment process at the end of this subsection.

Denote $v^n \approx v(t_n)$ to be the semi-discretization of $v$ in time, and also $v_h \approx v(\mathbf{u})$ to be the semi-discretization of $v$ in space. We write the space-time full discretization of $v$ as $v_h^n \approx v(\mathbf{u}^n,t_n)$.
The reformulated equation system \eqref{E:energy sav} can be discretized by the rational basis functions into the following form
\begin{align}\label{E:gBo sav dx}
\begin{cases}
\mathbf{u}_t=-\mathbf{P^{-1}F^{-1}S_1FP}\left(\mathbf{-P^{-1}F^{-1}HS_1FP u}+ \dfrac{1}{m}\frac{\mathbf{u^m} v_h}{\sqrt{\langle \mathbf{u^m},\mathbf{u} \rangle_h}+C_0} \right) \defeq f(\mathbf{u},v_h),\\
(v_h)_t=\frac{m+1}{2\sqrt{\langle \mathbf{u^m} ,\mathbf{u} \rangle_h+C_0}} \langle \mathbf{u^m} ,\mathbf{u}_t \rangle_h \defeq g(\mathbf{u},v_h),
\end{cases}
\end{align} 
with the initial conditions
$$\mathbf{u}^0=u(\mathbf{x},0), \quad \quad v_h^0=\sqrt{\langle (\mathbf{u}^0)^{\mathbf{m}},\mathbf{u}^0 \rangle_h+C_0}.$$
The fully discrete modified energy is defined as follows
\begin{align}\label{D:m energy sav}
E_h^n=\frac{1}{2} \langle \mathbf{P^{-1}F^{-1} HS_1 FP} \mathbf{u}^n, \mathbf{u}^n \rangle_h -\frac{1}{m(m+1)} \left((v_h^n)^2-C_0 \right).
\end{align}
Next, we prove the following theorem for the high order energy-conservative schemes.
\begin{theorem}\label{T:SRK sav}
The $s$-stage symplectic Runge-Kutta method, which satisfies \eqref{E:SRK ab}, 
conserves the $L^1$-type integral \eqref{E:momentum}, mass \eqref{E:mass} and modified energy \eqref{E:energy sav} in the discrete time flow for the reformulated gBO equation system \eqref{E:gBO sav}, i.e.,
\begin{align}\label{E:I M E c}
I^{n+1}=I^n, \quad M^{n+1}=M^{n}, \quad \mbox{and} \quad E^{n+1}=E^{n}.
\end{align}
Furthermore, the symplectic Runge-Kutta method preserves the discrete energy \eqref{D:m energy sav} for the spatial semi-discretized system \eqref{E:gBo sav dx}, i.e.,
\begin{align}\label{E:DE conservation}
E_h^n=E_h^{n-1}=\cdots=E_h^0.
\end{align}
\end{theorem}

\begin{proof}
The proof for the conservation of the temporal semi-discretized $L^1$-type integral, mass and energy \eqref{E:I M E c} is standard, e.g., see \cite{Cooper1987}, \cite{Sanz1988}, \cite{LZQS2019} and our earlier paper \cite{Yang2021}. 

For the proof of the discrete energy conservation \eqref{E:DE conservation}, substituting the inner product into the discrete sense, straightforward calculations yield
\begin{align*}
E^{n+1}_h&=\frac{1}{2} \langle \mathbf{P^{-1}F^{-1} HS_1 FP} \mathbf{u}^{n+1}, \mathbf{u}^{n+1} \rangle_h -\frac{1}{m(m+1)} \left((v_h^{n+1})^2-C_0 \right) \\
&=\frac{1}{2} \langle \mathbf{P^{-1}F^{-1} HS_1 FP} (\mathbf{u}^n+\tau \sum_{i=1}^{s} b_if(\mathbf{U}_i,V_i)), \mathbf{u}^n+\tau \sum_{i=1}^{s} b_if(\mathbf{U}_i,V_i) \rangle_h\\
&-\frac{1}{m(m+1)} \left((v_h^{n}+\tau \sum_{i=1}^{s} b_i g(\mathbf{U}_i,V_i))^2-C_0 \right)\\
&=E_h^n+ \tau \sum_{i=1}^s b_i \left( \langle \mathbf{P^{-1}F^{-1} HS_1 FP} f(\mathbf{U}_i,V_i),\mathbf{U}_i \rangle_h -\frac{2}{m(m+1)} V_i g(\mathbf{U}_i,V_i) \right)\\
&+\tau^2 \sum_{i,j=1}^s (b_ia_{ij}+b_ja_{ji}-b_ib_j) \left(\langle f(\mathbf{U}_i,V_i),f(\mathbf{U}_j,V_j) \rangle_h +g(\mathbf{U}_i,V_i)g(\mathbf{U}_j,V_j) \right)\\
&=E^n_h
\end{align*}
by using the relation \eqref{Eh identity}, \eqref{E:gBo sav dx} and \eqref{E:SRK ab}, where $V_i$ is defined as the intermediate value of $v_h^n$ for \eqref{E:gBo sav dx} similar to \eqref{Ui} and \eqref{u n+1}.

\end{proof}

\begin{remark}
Comparing with the proof in \cite[Theorem 3.1]{Yang2021}, the symplectic Runge-Kutta method can preserve all the three quantities in the discrete time flow for the reformulated system \eqref{E:gBO sav}. However, due to the limitation of the spatial discretization, only the discrete energy will be preserved in the fully discrete sense. The conservation of the discrete momentum for the gKdV equations in \cite{Yang2021} is obtained by using the circulant and anti-symmetric property of the first order differential matrix from the Fourier pseudo-spectral discretization. We note that the rational basis functions here do not possess this property. 
\end{remark}

The adjustment process for the constant $C_0$ from \cite{Yang2021} can be adapted here.  Suppose at $t=t_n$, the term  $\int (u^n)^{m+1} dx +C_0<Tol$, where $Tol$ is a given positive number (e.g., $Tol=5$). Then, we choose another constant $\tilde{C_0}$ such that $\int (u^n)^{m+1} dx +\tilde{C_0}> Tol$. For example, we can take $\tilde{C_0}=10-\int (u^n)^{m+1} dx$, which leads to our new $\tilde{v}^n \approx \sqrt{10}$. Then, by using $E[u^n,v^n]=E[u^n,\tilde{v}^n]$ from \eqref{E:energy sav}, we have our new $\tilde{v}^n$
\begin{align}\label{E:v new}
\tilde{v}^n=\sqrt{(v^n)^2+\tilde{C_0}-C_0}.
\end{align}
Finally, we substitute the $v^n$ and $C_0$ in \eqref{E:gBO sav} with $\tilde{v}^n$ and $\tilde{C_0}$, and then, continue with the time evolution for $t=t_{n+1},t_{n+2}, \cdots$.
 
\begin{remark}
Note that $v^2=\int u^{n+1} dx +C_0$ holds only at the collocation points $t=t_n+\tau c_i$ for each $i=1,2,\cdots,s$ in $t\in [t_n,t_{n+1}]$. However, the constant $c_i$ may not necessarily be equal to $0$ or $1$, e.g., see Table \ref{T:RK Butcher}. This means $v^2=\int u^{n+1} dx +C_0$ does not hold at $t_n$ in the discrete time flow.
Therefore, the new $\tilde{v}^n$ can only be evaluated by \eqref{E:v new} to keep the discrete energy \eqref{E:energy sav} invariant.
\end{remark}

\section{Numerical results}
In this section, we list our numerical examples for the proposed schemes. Before discussing the examples, we mention that a type of a fixed point iteration solver from 
\cite{CWJ2021} and \cite{Yang2021} can be easily adapted here for solving the resulting nonlinear system from the IRK methods. The total computational cost is on the order of $ \mathcal{O}(N\log(N))$ from FFT.

We denote by the IRK2-MC and IRK4-MC the mass conservative schemes for solving \eqref{E:gBO SMC} by using the 1st and 2nd stage RK methods with Gauss-Legendre collocation points from Table \ref{T:RK Butcher}. We also denote by the IRK2-EC and IRK4-EC the energy conservative schemes for solving the reformulated system \eqref{E:gBo sav dx}. As a comparison, we use the commonly used 2nd order non-conservative semi-implicit Leap-Frog scheme as follows
\begin{align}\label{E:Leap-Frog}
\frac{\mathbf{u}^{n+1}-\mathbf{u}^{n-1}}{2 \tau} = -\mathbf{P^{-1}F^{-1}S_1FP}\left(\mathbf{P^{-1}F^{-1}HS_1FP} \left(\frac{\mathbf{u}^{n+1}+\mathbf{u}^{n-1}}{2} \right) + \frac{1}{m}(\mathbf{u}^{n})^{\mathbf{m}} \right),
\end{align}
denoted as Leap-Frog.
We track the following quantities at $t=t_n$ to check the accuracy:
\begin{align}
&\mathcal{E}^n=\|u_{exact}^n-\mathbf{u}^n\|_{\infty}; \label{E:error inf} \\
&\mathcal{E}_I^n=\max_{l<n}|I_h^l-I_h^0|; \label{E:error I} \\
&\mathcal{E}^n_M=\max_{l<n}|M_h^l-M_h^0|; \label{E:error M} \\
&\mathcal{E}^n_E=\max_{l<n}|E_h^l-E_h^0|. \label{E:error E}
\end{align}

When the SAV approach is not applied (IRK2-MC , IRK4-MC and Leap-Frog), the discrete energy $E_h^n$ is computed from \eqref{D:energy st}; and when the SAV approach is applied, the discrete energy $E_h^n$ is computed from the modified version \eqref{D:m energy sav}.
We mention here that it is easy to see the equivalence between the Crank-Nicholson scheme \eqref{CN MC} and the IRK2-MC scheme. The energy-conservative Crank-Nicholson scheme, which \eqref{CN EC} considers reformulating the potential, also shares the same idea as in the IRK2-EC scheme. We omit the numerical result from the Crank-Nicholson methods for the purpose of conciseness, though compared with the IRK2-EC scheme, the energy-conservative Crank-Nicholson scheme usually requires less iterations in the fixed point iteration process, since the SAV approach introduces an additional scalar variable.

Now, we are ready to illustrate examples for our numerical simulations.

\medskip
{\flushleft \it \underline{Example $1$}}. Our first example considers the soliton solution for the BO ($m=2$) equation, $u(x,t)=\frac{4c}{1+c^2(x-x_0-ct)^2}$. These type of solutions come from the smooth, positive, decaying at infinity solitary wave solution to the profile equation
\begin{align}\label{E:GS c}
\mathcal{H}Q_x+cQ-\frac{1}{m}Q^m=0,
\end{align}
where $c$ is a constant that indicates the speed of the traveling waves as well as the magnititude. The solutions are expected to travel to the right as the solitons, for example from the spectral stability result in \cite{KM2009}, \cite{ABH} and the inverse scattering theory \cite{FA_83}, \cite{KLM_98}.

In our numerical simulations, we take $\alpha=25$ with $N=1024$. We take the traveling speed $c=2$ and starting point at $x_0=-20$. The time step $\tau$ is taken to be $\tau=\frac{1}{20}$ for all the four IRK type methods, and $\tau=\frac{1}{40}$ for the Leap-Frog scheme \eqref{E:Leap-Frog}, since taking the $\tau=\frac{1}{20}$ will lead to the numerical instability in our numerical computations for the Leap-Frog scheme. We stop our numerical simulation at $T=20$.

\begin{figure}
\includegraphics[width=0.45\textwidth]{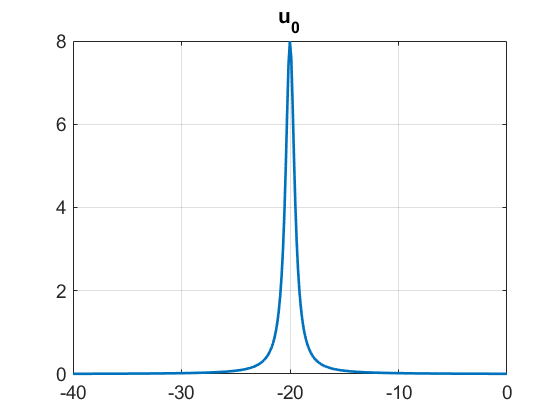}
\includegraphics[width=0.45\textwidth]{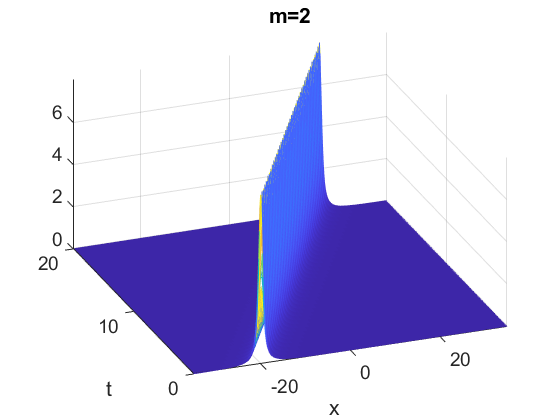}
\caption{\label{F:EX1 profile} The solution profile for Example 1 from the IRK4-EC. Left: $u_0$. Right: $u(x,t)$.}
\end{figure}

Figure \ref{F:EX1 profile} shows the solution profile obtained from the IRK4-EC scheme. The left subplot is the initial condition $u_0$, the right subplot is the time evolution. One can see that the solution travels in the solitary wave manner, which is previously observed in \cite{JW92} and \cite{RWY2021} and also as expected. 

\begin{figure}
\includegraphics[width=0.45\textwidth]{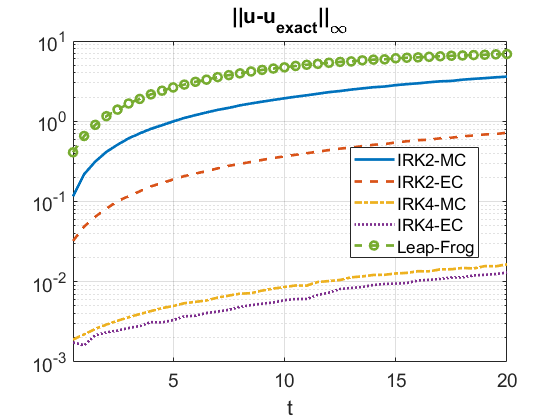}
\includegraphics[width=0.45\textwidth]{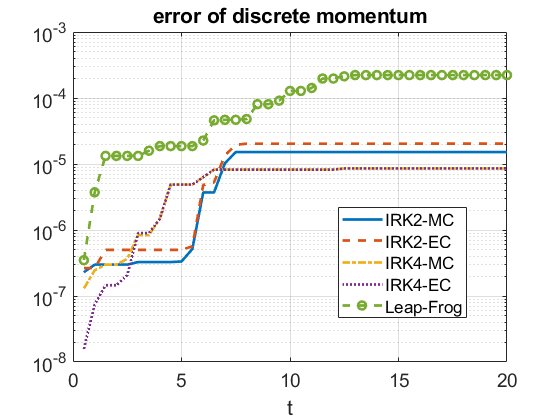}
\includegraphics[width=0.45\textwidth]{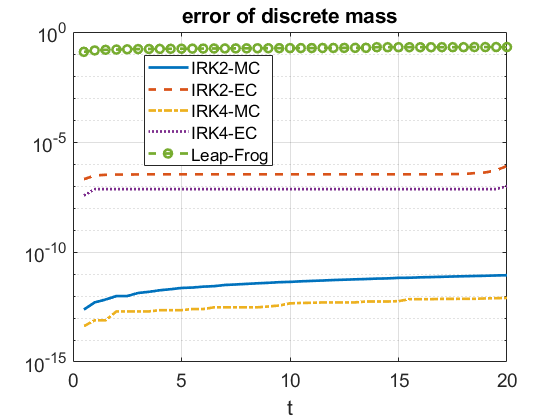}
\includegraphics[width=0.45\textwidth]{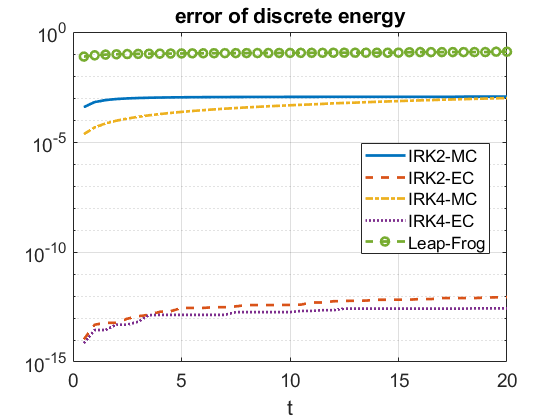}
\caption{\label{F:EX1 data} The errors in Example 1 by different time integrators: IRK2-MC (solid blue); IRK2-EC(dash red); IRK4-MC (dash dot orange); IRK4-EC (dot purple); Leap-Frog (circle green). Top left: $\|u-u_{exact}\|_{\infty}$. Top right: discrete momentum error. Bottom left: discrete mass error. Bottom right: discrete energy error. }
\end{figure}

Figure \ref{F:EX1 data} tracks the results obtained from the different time integrators. 
The top left subplot in Figure \ref{F:EX1 data} shows $\|\mathbf{u}^n-u_{exact}\|_{\infty}$ with respect to time, where $u_{exact}=u(\mathbf{x},t_n)$ is the exact solution. One can see that the Leap-Frog scheme \eqref{E:Leap-Frog} has the largest error (see the green circle line). Moreover, the 4th order schemes (IRK4-MC and IRK4-EC) own the better accuracy than the second order schemes (IRK2-MC, IRK2-EC and Leap-Frog), which is as expected, since they have higher order temporal accuracy. Furthermore, we observe that the energy-conservative schemes (dash red line for IRK2-EC and dot purple line for IRK4-EC) perform better than the 
mass-conservative schemes (blue solid line for IRK2-MC and orange dash-dot line for IRK4-EC). 
Thus, the energy-conservative schemes are recommended for future studies.

The top right subplot in Figure \ref{F:EX1 data} tracks the error of the discrete $L^1$-type integral \eqref{E:error I} at different times. One can see that the more accurate the time integrators are, the better preservation of this quantity will be, though they are not conserved exactly.

The bottom two subplots in Figure \ref{F:EX1 data} track the error of discrete mass \eqref{E:error M} and energy \eqref{E:error E}, respectively. The mass-conservative schemes (IRK2-MC and IRK4-MC) keep the error of discrete mass below at the level of $10^{-12}$, which is the tolerance of the fixed iteration in solving the resulting nonlinear system from the implicit Runge-Kutta method. On the other hand, the discrete mass error for the other types of schemes is relatively large, especially the Leap-Frog scheme (thus, refered to as least accurate). 

The bottom right subplot tracks the error of the discrete energy. It shows that the energy-conservative schemes (IRK2-EC and IRK4-EC) keep the error of discrete mass below the level of $10^{-12}$. This justifies the validity of our schemes. Similarly, for the other time integrators, the Leap-Frog scheme performs the worst even with a smaller time step $\tau$,; the other two mass-conservative schemes keep the error of discrete energy around the level of $10^{-4}$.

We also list the $L^{\infty}$ error $\mathcal{E}^n$ at $t=T$ with different time step $\tau$ for these five time integrators in Table \ref{T:EX1 error}. One can see that the Leap-Frog, IRK2-MC and IRK2-EC decrease with the ratio around $4$, which are as expected, since they are of the second order schemes. On the other hand, the ratio of the 4th order methods IRK4-MC and IRK4-EC is around $2^4$, since they are 4th order methods. When the time step $\tau$ is small, the decay rate is slightly below $16$ (see the last row in Table \ref{T:EX1 error}). This is probably because the temporal error becomes comparable to the spatial discretization error, and consequently, affects the ratio.

\begin{small}
\begin{table}
\begin{tabular}{|c|c|c|c|c|c|c|c|c|c|c|}
\hline
 & \multicolumn{2}{|c|}{Leap-Frog} &\multicolumn{2}{|c|}{IRK2-MC} &  \multicolumn{2}{|c|}{IRK2-EC} &\multicolumn{2}{|c|}{IRK4-MC} & \multicolumn{2}{|c|}{IRK4-EC}\\
\hline
  $\tau$& error & rate &  error  & rate & error  & rate & error  & rate & error  & rate\\
  \hline
 $\frac{1}{10}$& NA & NA & $7.05$ & NA  &$3.45$ & NA &$0.35$ & NA &$0.61$ & NA  \\
 \hline
 $\frac{1}{20}$& NA & NA & $3.60$& $1.96$ &$0.72$ &  $4.80$ &$1.6e-2$& $21.6$&$1.3e-2$ &  $46.8$ \\
 \hline
 $\frac{1}{40}$& $6.88$ & NA &$1.02$ & $3.53$ &$0.18$ & $4.05$  &$8.9e-4$ & $18.6$ &$5.6e-4$ & $23.2$ \\
 \hline
 $\frac{1}{80}$& $1.92$ & $3.59$ &$0.26$ & $3.93$ &$0.44$ & $4.03$  &$6.7e-5$ & $13.2$ &$5.1e-5$ & $10.9$  \\
 \hline
\end{tabular}
\caption{\label{T:EX1 error} The convergence rates of Leap-Frog, IRK2-MC, IRK2-EC IRK4-MC and IRK4-EC in Example 1.}
\end{table}
\end{small}

\medskip
{\flushleft \it \underline{Example $2$}}. We next consider the scattering solution for the BO equation with the initial condition $u_0=- 2 \sech^2(x)$. Its KdV version has been studied for questions on dispersion limit, see, e.g. \cite{GK2012} and \cite{KP2015}. Here, we expect that a similar solution behavior may happen, since the BO equation only changes the dispersion term $u_{xxx}$ from the KdV equation to $\mathcal{H} u_{xx}$ (less amount of dispersion if viewed on the Fourier frequency side). Note that a negative value for $\int u^3 dx$ may occur, and thus, the $C_0$ adjustment process in \eqref{E:v new} will make $v(t)$ stay positive, and will keep the algorithm applicable for all time. 
The exact solution is not explicitly given, since due to the negative sign in the initial condition and coefficients chosen. 

\begin{figure}
\includegraphics[width=0.7\textwidth]{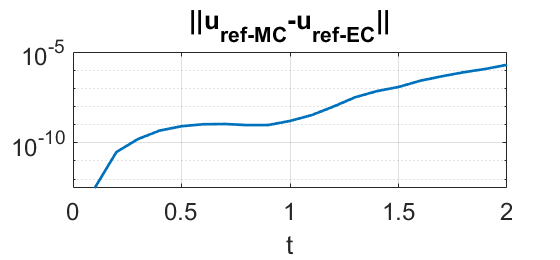}
\caption{\label{F:EX2 uref} The difference of the reference solution obtained by the IRK4-EC scheme and IRK4-MC scheme. One can see the difference (mainly caused by the spatial discretization) keeps increasing to the level $10^{-6}$.}
\end{figure}

In this example, we still take the $N=1024$ and $\alpha=25$ for the spatial discretization. The time step $\tau=\frac{1}{400}$ ($\tau=\frac{1}{800}$ for the Leap-Frog) and the stopping time $T=2$. We compute the reference solution $u_{\mathrm{ref}}$ by both IRK4-MC and IRK4-EC methods independently with an ultimately small time step ($\tau=1/6400$), denoted as $u_{\mathrm{ref-MC}}$ and $u_{\mathrm{ref-EC}}$, respectively. Since we intend to track the convergence rate with respect to time, to minimize the influence from the spatial discretization error, we use the $u_{\mathrm{ref-MC}}$ to compute the $L^{\infty}$ error $\| \mathbf{u}^n-u_{\mathrm{ref}} \|_{\infty} $ when the $\mathbf{u}^n$ is obtained by the mass-conservative schemes (IRK2-MC and IRK4-MC), and use the $u_{\mathrm{ref-EC}}$ to compute the $L^{\infty}$ error $\|\mathbf{u}^n-u_{\mathrm{ref}}\|_{\infty}$ when $\mathbf{u}^n$ is obtained by the energy-conservative schemes (IRK2-EC and IRK4-EC) and the Leap-Frog scheme. The spatial discretization error accumulates as the time evolves, see Figure \ref{F:EX2 uref}. From Figure \ref{F:EX2 uref}, the difference increases to the level of $10^{-6}$ between these two solutions by the time we terminate the simulation. 

\begin{figure}
\includegraphics[width=0.45\textwidth]{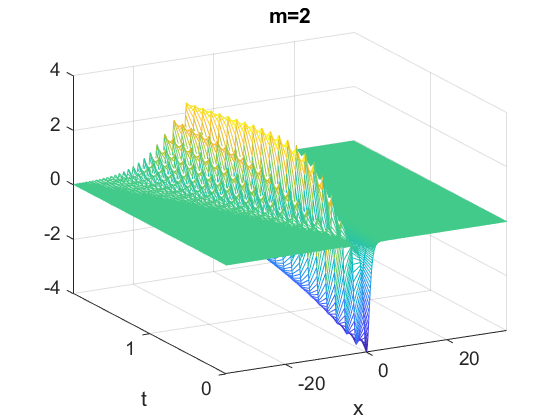}
\includegraphics[width=0.45\textwidth]{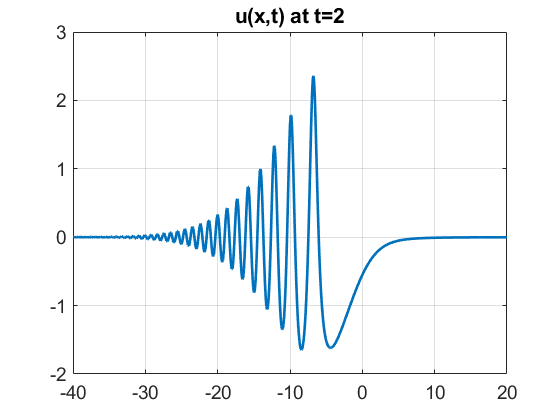}
\caption{\label{F:EX2 profile} The solution profile in Example 2 from IRK4-SAV. Left: $u(x,t)$. Right: $u(x,t)$ at $t=2$.}
\end{figure}

Figure \ref{F:EX2 profile} shows the solution profile obtained from the IRK4-EC method. The left subplot shows the solution profile at different times $t$. The right plot shows the solution at the terminal time $t=2$. We can see that the solution radiates to the right with fast oscillations. On the other hand, compared with the similar type of solutions to the KdV case (e.g., in \cite{KS2015} and \cite{Yang2021}), the frequency is smaller. This indicates that the lower order dispersion ($\mathcal{H}\partial_{xx}$ compared with $\partial_{xxx}$) generates slower oscillations.

\begin{figure}
\includegraphics[width=0.45\textwidth]{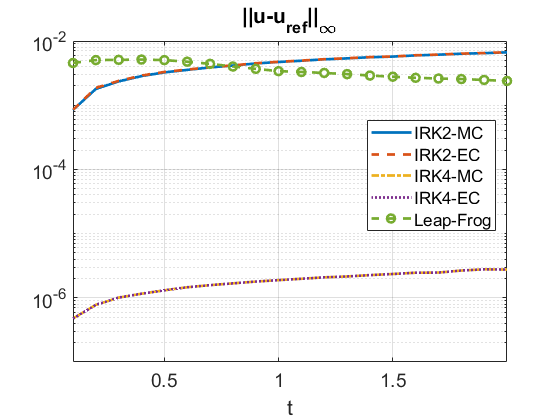}
\includegraphics[width=0.45\textwidth]{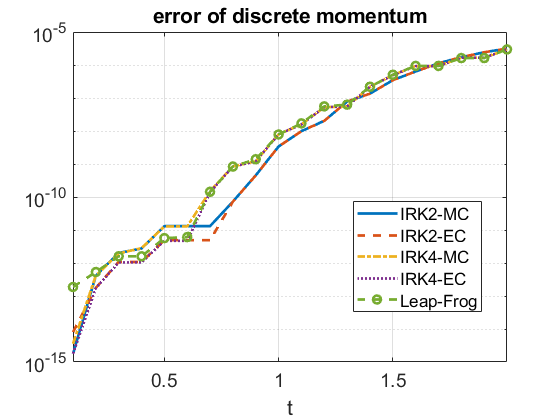}
\includegraphics[width=0.45\textwidth]{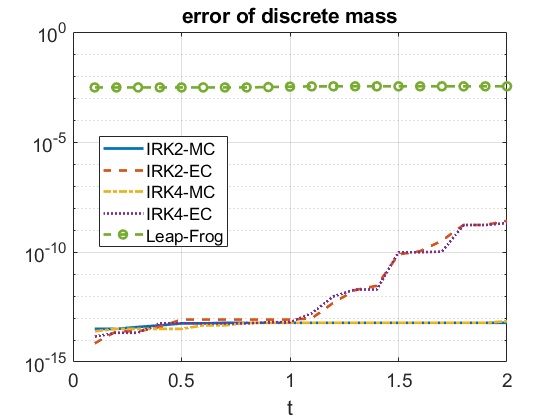}
\includegraphics[width=0.45\textwidth]{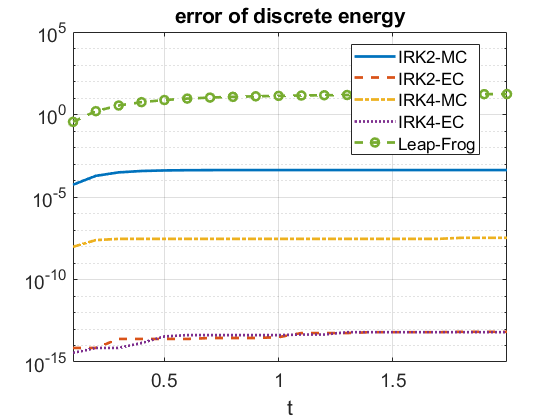}
\caption{\label{F:EX2 data} The errors in Example 2 by different time integrators: IRK2-MC (solid blue); IRK2-EC(dash red); IRK4-MC (dash dot orange); IRK4-EC (dot purple); Leap-Frog (circle green). Top left: $\| \mathbf{u}^n-u_{\mathrm{ref}}\|_{\infty}$. Top right: discrete momentum error. Bottom left: discrete mass error. Bottom right: discrete energy error.}
\end{figure}

Figure \ref{F:EX2 data} tracks the $L^{\infty}$-error, error of discrete $L^1$-type integral, mass and energy with respect to time. 
One can see that the results are similar to the previous example, and also agree with our analysis in Section \ref{S:SD} and \ref{S:TD}.

Table \ref{T:EX2 error} shows the $L^{\infty}$ error at $t=T$ with respect to the different time step $\tau$. The decay rate is on the 2nd order for the 2nd order schemes (Leap-Frog, IRK2-MC and IRK2-EC), and on the 4th order for the 4th order schemes (IRK4-MC and IRK4-EC). Surprisingly, the IRK types of schemes (IRK2-MC and IRK2-EC, IRK4-MC and IRK4-EC) generate almost the same error (up to the decimals that we report) from the different reference solutions ($u_{\mathrm{ref-MC}}$ and $u_{\mathrm{ref-EC}}$). This implies that some possible cancellations may occur between the spatial discretization errors. 
 
\begin{small}
\begin{table}
\begin{tabular}{|c|c|c|c|c|c|c|c|c|c|c|}
\hline
 & \multicolumn{2}{|c|}{Leap-Frog} &\multicolumn{2}{|c|}{IRK2-MC} &  \multicolumn{2}{|c|}{IRK2-EC} &\multicolumn{2}{|c|}{IRK4-MC} & \multicolumn{2}{|c|}{IRK4-EC}\\
\hline
  $\tau$& error & rate &  error  & rate & error  & rate & error  & rate & error  & rate\\
  \hline
 $\frac{1}{200}$& $0.104$ & NA & $0.027$ & NA  &$0.027$ & NA &$4.4e-5$ & NA &$4.4e-5$ & NA  \\
 \hline
 $\frac{1}{400}$& $0.026$ & $4.00$ & $6.7e-3$& $3.99$ &$6.7e-3$ &  $3.99$ &$2.8e-6$& $15.8$&$2.8e-6$ &  $15.8$ \\
 \hline
 $\frac{1}{800}$& $6.5e-3$ & $3.99$ &$1.7e-3$ & $4.02$ &$1.7e-3$ & $4.02$  &$1.7e-7$ & $16.0$ &$1.7e-7$ & $16.0$ \\
 \hline
 $\frac{1}{1600}$& $1.6e-3$ & $4.00$ &$4.2e-4$ & $3.99$ &$4.2e-4$ & $3.99$  &$1.1e-8$ & $16.0$ &$1.1e-8$ & $16.0$  \\
 \hline
\end{tabular}
\caption{\label{T:EX2 error} The convergence rates of Leap-Frog, IRK2-MC, IRK2-EC IRK4-MC and IRK4-EC in Example 2.}
\end{table}
\end{small} 

\medskip
{\flushleft \it \underline{Example $3$}}. Our final example considers the mBO ($m=3$) and the gBO ($m=4$) cases. We take the initial condition $u_0=0.99Q$, where $Q$ is the soliton solution from \eqref{E:GS c} with $c=1$. In these cases, while there is no explict form for $Q$, the profile of $Q$ can be obtained numerically, e.g., by the Petviashvili iteration from \cite{PS2004}, \cite{RWY2021}, and its convergence analysis in \cite{LY2007}, \cite{OSSS2016} and \cite{LP2019}.
From \cite{LFA2014}, when $u_0 =0.99Q$, which indicates that the solution is below the mass-energy threshold, the solution is proven to exist globally in time.
Recent numerical study in \cite{RWY2021} shows that the solution blows up when $u_0=1.01Q$. In this paper, we consider the globally existing solutions, and thus, we take $u_0=0.99Q$ in our example. We take $N=1024$, $\alpha=25$, $\tau=0.02$ ($\tau=0.01$ for the Leap-Frog scheme due to the stability issue) in our simulation. We run until $T=10$ for $m=3$, and $T=5$ for $m=4$.

\begin{figure}
\includegraphics[width=0.45\textwidth]{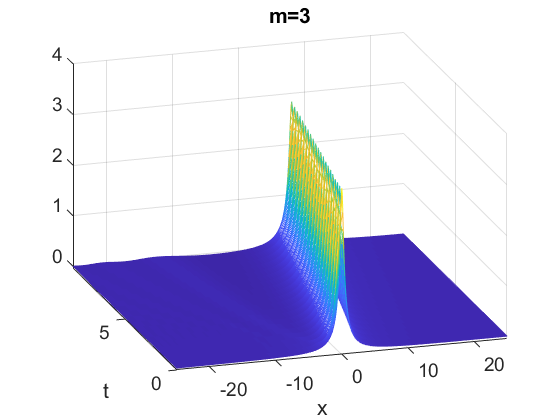}
\includegraphics[width=0.45\textwidth]{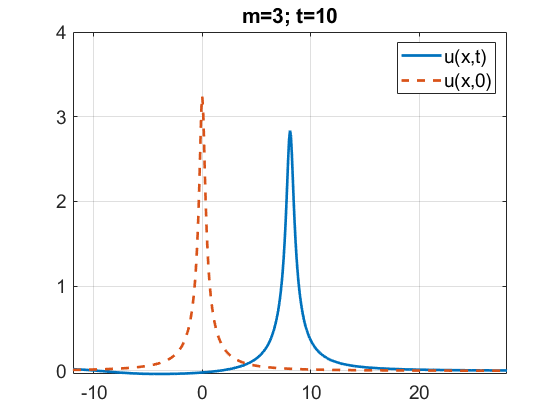}
\caption{\label{F:EX3 profile m3} The solution profile in Example 3 from IRK4-EC with $m=3$. Left: $u(x,t)$. Right: $u(x,t)$ at $t=10$.}
\end{figure}

\begin{figure}
\includegraphics[width=0.45\textwidth]{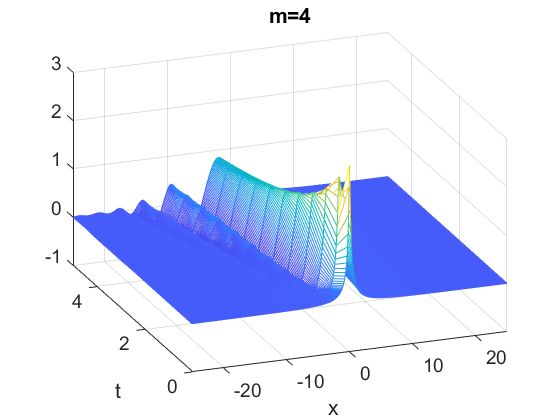}
\includegraphics[width=0.45\textwidth]{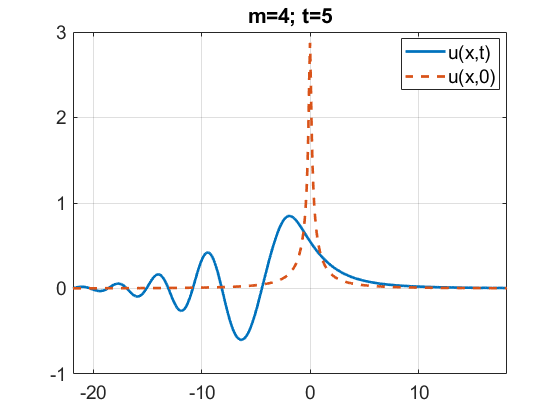}
\caption{\label{F:EX3 profile m4} The solution profile in Example 3 from IRK4-EC with $m=4$. Left: $u(x,t)$. Right: $u(x,t)$ at $t=5$.}
\end{figure}

Figures \ref{F:EX3 profile m3} and \ref{F:EX3 profile m4} show the solution profiles obtained from the scheme IRK4-EC for $m=3$ and $m=4$, respectively. The left subplots show the solution profiles $u(x,t)$ at different times. The right subplots show the solution at the final time $T$ (blue solid line) and their comparison with the initial condition $u_0$ (red dash line). For $m=3$, the solution travels to the right with some radiation parts scattering to the left. However, for the $m=4$ case, the solution completely radiates to the left. This agrees with the results in our earlier paper \cite{RWY2021}.

\begin{figure}
\includegraphics[width=0.32\textwidth]{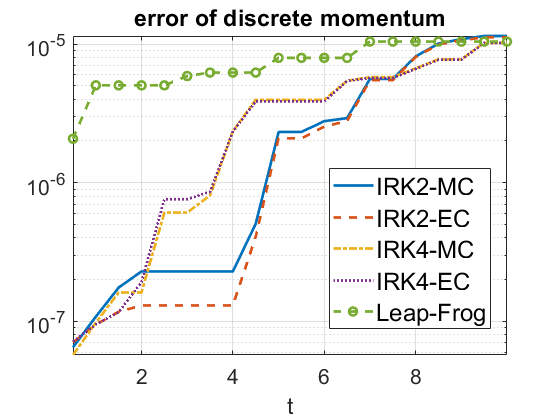}
\includegraphics[width=0.32\textwidth]{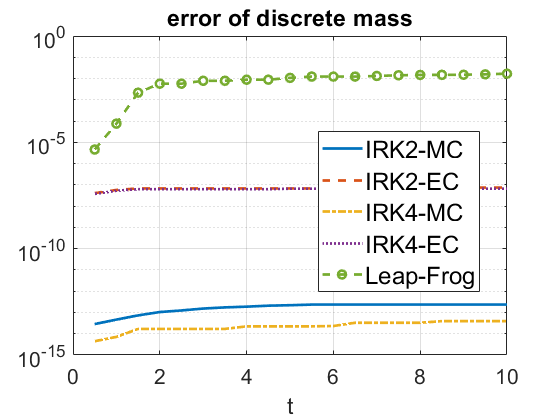}
\includegraphics[width=0.32\textwidth]{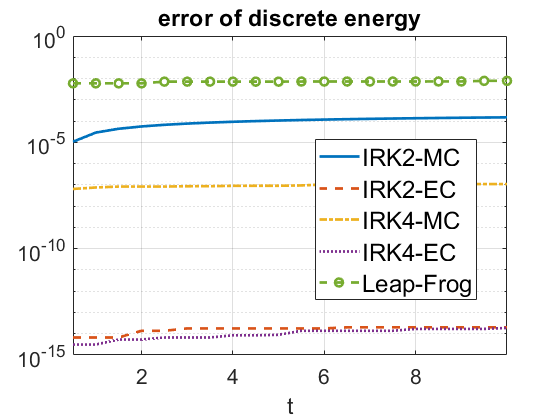}
\caption{\label{F:EX3 data m3} The errors in Example 3 ($m=3$) by different time integrators: IRK2-MC (solid blue); IRK2-EC(dash red); IRK4-MC (dash dot orange); IRK4-EC (dot purple); Leap-Frog (circle green). Left: discrete momentum error. Middle: discrete mass error. Right: discrete energy error.}
\end{figure}

\begin{figure}
\includegraphics[width=0.32\textwidth]{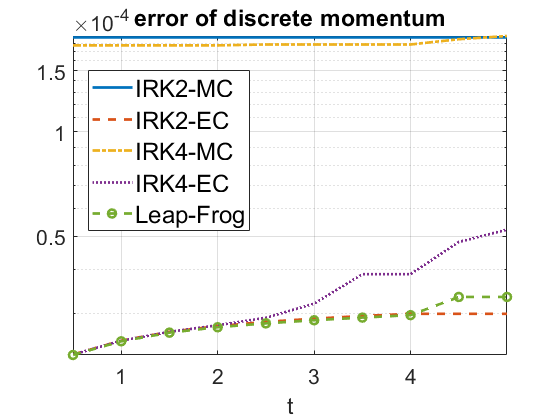}
\includegraphics[width=0.32\textwidth]{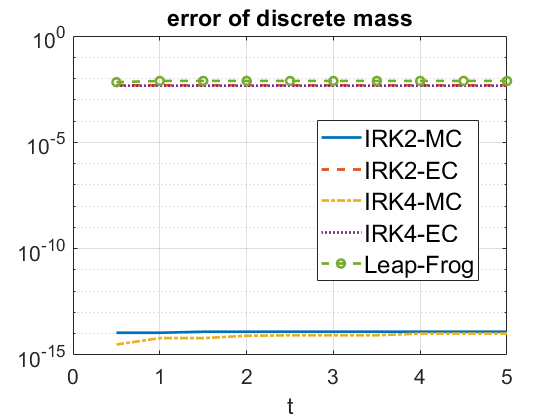}
\includegraphics[width=0.32\textwidth]{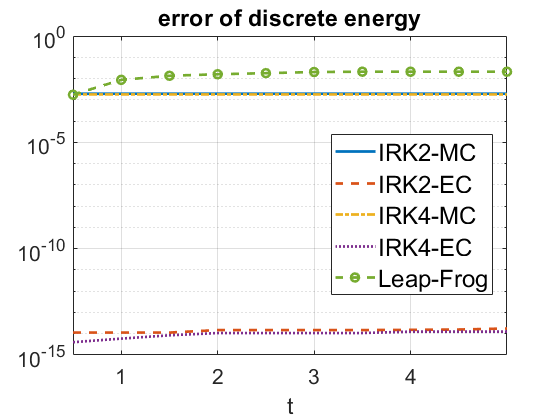}
\caption{\label{F:EX3 data m4} The errors in Example 3 ($m=4$) by different time integrators: IRK2-MC (solid blue); IRK2-EC(dash red); IRK4-MC (dash dot orange); IRK4-EC (dot purple); Leap-Frog (circle green). Left: discrete momentum error. Middle: discrete mass error. Right: discrete energy error.}
\end{figure}

Figure \ref{F:EX3 data m3} and \ref{F:EX3 data m4} track the error of discrete $L^1$-type integral (left subplot), mass (middle subplot) and energy (right subplot) at different times for $m=3$ and $m=4$, respectively. 
Again, the discrete mass or energy can be preserved by choosing the mass-conservative scheme or energy-conservative scheme, respectively, which agrees with the analysis in Sections \ref{S:SD} and \ref{S:TD}.

\section{Conclusion and other discussion}\label{S:conclusion}
We propose two kinds of pseudo-spectral spatial discretization for the generalized Benjamin-Ono equation, one is mass-conservative, and the other one is energy-conservative. 
Combined with the conservative time discretization, such as the symplectic Runge-Kutta method, arbitrarily high order mass-conservative or energy-conservative numerical schemes can be constructed.
In fact, the symplectic Runge-Kutta method with the scalar auxiliary variable reformulation preserve all three invariants \eqref{E:momentum}-\eqref{E:energy} in the discrete time flow. However, the spatial discretization restricts us constructing the scheme to preserve more than one invariant quantity in the fully discrete sense.
Numerical results show that the energy-conservative schemes possess superior accuracy compared to the mass-conservative schemes, nevertheless, both of these schemes are more accurate than the non-conservative (Leap-Frog) scheme.

This strategy can be extend to construct the the conservative schemes for the gKdV equations, i.e.,
$$u_t=-u_{xxx}-\frac{1}{m}(u^m)_x;$$
and also the mass-energy conservative (structure-preserving) schemes for the NLS equations, i.e., 
$$u_t=i\left(u_{xx}+|u|^{m-1}u \right).$$ 
For example, the corresponding spatial semi-discretized form of the equation \eqref{E:gBo sav dx} for the gKdV equation will be  
\begin{align*}
\begin{cases}
\mathbf{u}_t=-\mathbf{P^{-1}F^{-1}S_1FP}\left(\mathbf{-P^{-1}F^{-1}S_2FP u}+ \dfrac{1}{m}\frac{\mathbf{u^m} v_h}{\sqrt{\langle \mathbf{u^m},\mathbf{u} \rangle_h}+C_0} \right),\\
(v_h)_t=\frac{m+1}{2\sqrt{\langle \mathbf{u^m} ,\mathbf{u} \rangle_h+C_0}} \langle \mathbf{u^m} ,\mathbf{u}_t \rangle_h.
\end{cases}
\end{align*} 
Similarly, the semi-discretized form for the NLS equation will be
\begin{align*}
\begin{cases}
\mathbf{u}_t=i\left(\mathbf{-P^{-1}F^{-1}S_2FP u}+ \frac{\mathbf{|u|^{m-1}u} v_h}{\sqrt{\langle \mathbf{|u|^{m-1}u},\mathbf{u} \rangle_h}+C_0} \right),\\
(v_h)_t=\frac{m}{2\sqrt{\langle \mathbf{|u|^{m-1}u} ,\mathbf{u} \rangle_h+C_0}} \Re \left( \langle \mathbf{|u|^{m-1}u} ,\mathbf{u}_t \rangle_h \right).
\end{cases}
\end{align*} 
The straightforward adaption of proof in Theorem \ref{T:SRK sav} will show the energy-conservative result for the gKdV equations, and the structure-preserving result for the NLS equations. These results can also be easily extended to higher dimensions (e.g., Zakharov–Kuznetsov equation or the $d$-dimensional NLS equation) by applying the tensor product in the spatial discretization. We omit the proof and numerical examples here for conciseness.

In summary, by applying the rational basis functions, the above illustrated conservative schemes will increase the computational efficiency significantly, especially, in tracking the solution's long time behavior, or the slow decaying solutions (e.g., $u_0=\frac{1}{1+x^2}$), since far less number of nodes are needed compared with the traditional domain truncation approaches.

\bibliographystyle{abbrv}
\bibliography{ref_BOc}

\begin{thebibliography}{10}

\bibitem{ABFS1989}
L.~Abdelouhab, J.~Bona, M.~Felland, and J.-C. Saut.
\newblock Nonlocal models for nonlinear, dispersive waves.
\newblock {\em Physica D: Nonlinear Phenomena}, 40(3):360--392, 1989.

\bibitem{AC1991}
M.~J. Ablowitz, M.~Ablowitz, P.~Clarkson, and P.~A. Clarkson.
\newblock {\em Solitons, nonlinear evolution equations and inverse scattering},
  volume 149.
\newblock Cambridge university press, 1991.

\bibitem{AFA_83}
M.~J. Ablowitz, A.~S. Fokas, and R.~L. Anderson.
\newblock The direct linearizing transform and the {B}enjamin-{O}no equation.
\newblock {\em Phys. Lett. A}, 93(8):375--378, 1983.

\bibitem{ABH}
J.~Albert, J.~L. Bona, and D.~Henry.
\newblock Sufficient conditions for stability of solitary-wave solutions of
  model equations for long waves.
\newblock {\em Physica D}, 24:343--366, 1987.

\bibitem{benjamin_1967}
T.~B. Benjamin.
\newblock Internal waves of permanent form in fluids of great depth.
\newblock {\em Journal of Fluid Mechanics}, 29(3):559--592, 1967.

\bibitem{BBM1972}
T.~B. Benjamin, J.~L. Bona, and J.~J. Mahony.
\newblock Model equations for long waves in nonlinear dispersive systems.
\newblock {\em Philosophical Transactions of the Royal Society of London.
  Series A, Mathematical and Physical Sciences}, 272(1220):47--78, 1972.

\bibitem{BK_79}
T.~L. Bock and M.~D. Kruskal.
\newblock A two-parameter {M}iura transformation of the {B}enjamin-{O}no
  equation.
\newblock {\em Phys. Lett. A}, 74(3-4):173--176, 1979.

\bibitem{KB2000}
J.~Bona and H.~Kalisch.
\newblock {Models for internal waves in deep water}.
\newblock {\em Discrete \& Continuous Dynamical Systems-A}, 6(1):1--20, 2000.

\bibitem{BK04}
J.~Bona and H.~Kalisch.
\newblock Singularity formation in the generalized {B}enjamin-{O}no equation.
\newblock {\em Discrete and continuous dynamical systems}, 11:27--46, 2004.

\bibitem{BSS1987}
J.~L. Bona, P.~E. Souganidis, and W.~A. Strauss.
\newblock Stability and instability of solitary waves of {K}orteweg-de {V}ries
  type.
\newblock {\em Proceedings of the Royal Society of London. A. Mathematical and
  Physical Sciences}, 411(1841):395--412, 1987.

\bibitem{Bo87}
J.~P. Boyd.
\newblock Spectral methods using rational basis functions on an infinite
  interval.
\newblock {\em Journal of Computational Physics}, 69(1):112--142, 1987.

\bibitem{Boyd1990}
J.~P. Boyd.
\newblock The orthogonal rational functions of {H}iggins and {C}hristov and
  algebraically mapped {C}hebyshev polynomials.
\newblock {\em Journal of Approximation Theory}, 61(1):98--105, 1990.

\bibitem{BX11}
J.~P. Boyd and Z.~Xu.
\newblock Comparison of three spectral methods for the {B}enjamin-{O}no
  equation: {F}ourier pseudospectral, rational {C}hristov functions and
  {G}aussian radial basis functions.
\newblock {\em Wave Motion}, 48(8):702--706, 2011.

\bibitem{BX12}
J.~P. Boyd and Z.~Xu.
\newblock Numerical and perturbative computations of solitary waves of the
  {B}enjamin-{O}no equation with higher order nonlinearity using {C}hristov
  rational basis functions.
\newblock {\em J. Comput. Phys.}, 231(4):1216--1229, 2012.

\bibitem{BP2006}
N.~Burq and F.~Planchon.
\newblock Smoothing and dispersive estimates for 1{D} {S}chr\"odinger equations
  with {BV} coefficients and applications.
\newblock {\em J. Funct. Anal.}, 236(1):265--298, 2006.

\bibitem{BP2008}
N.~Burq and F.~Planchon.
\newblock On the well-posedness of the {B}enjamin-{O}no equation.
\newblock {\em Math. Ann.}, 340:497–--542, 2008.

\bibitem{Butcher1964}
J.~C. Butcher.
\newblock Implicit {R}unge-{K}utta processes.
\newblock {\em Mathematics of Computation}, 18(85):50--64, 1964.

\bibitem{Christov82}
C.~I. Christov.
\newblock A complete orthonormal system of functions in {$L^{2}(-\infty
  ,\,\infty )$} space.
\newblock {\em SIAM J. Appl. Math.}, 42(6):1337--1344, 1982.

\bibitem{Cooper1987}
G.~Cooper.
\newblock Stability of {R}unge-{K}utta methods for trajectory problems.
\newblock {\em IMA journal of numerical analysis}, 7(1):1--13, 1987.

\bibitem{CWJ2021}
J.~Cui, Y.~Wang, and C.~Jiang.
\newblock Arbitrarily high-order structure-preserving schemes for the
  {G}ross--{P}itaevskii equation with angular momentum rotation.
\newblock {\em Computer Physics Communications}, 261:107767, 2021.

\bibitem{davis_acrivos_1967}
R.~E. Davis and A.~Acrivos.
\newblock Solitary internal waves in deep water.
\newblock {\em Journal of Fluid Mechanics}, 29(3):593--607, 1967.

\bibitem{DM2002}
A.~Debussche and L.~{Di Menza}.
\newblock Numerical simulation of focusing stochastic nonlinear {S}chr\"odinger
  equations.
\newblock {\em Physica D: Nonlinear Phenomena}, 162(3):131--154, 2002.

\bibitem{DHKR16}
R.~Dutta, H.~Holden, U.~Koley, and N.~H. Risebro.
\newblock Convergence of finite difference schemes for the {B}enjamin--{O}no
  equation.
\newblock {\em Numerische Mathematik}, 134(2):249--274, Oct 2016.

\bibitem{DM2008}
T.~Duyckaerts and F.~Merle.
\newblock Dynamics of threshold solutions for energy-critical wave equation.
\newblock {\em Int. Math. Res. Pap. IMRP}, pages Art ID rpn002, 67, 2008.

\bibitem{LFA2014}
L.~G. Farah, F.~Linares, and A.~Pastor.
\newblock Global well-posedness for the $k$-dispersion generalized
  {B}enjamin-{O}no equation.
\newblock {\em Differential Integral Equations}, 27(7/8):601--612, 07 2014.

\bibitem{FA_83}
A.~S. Fokas and M.~J. Ablowitz.
\newblock The inverse scattering transform for the {B}enjamin-{O}no
  equation---a pivot to multidimensional problems.
\newblock {\em Stud. Appl. Math.}, 68(1):1--10, 1983.

\bibitem{Geng1993}
S.~Geng.
\newblock Construction of high order symplectic {R}unge-{K}utta methods.
\newblock {\em Journal of Computational Mathematics}, pages 250--260, 1993.

\bibitem{GWWC2017}
Y.~Gong, Q.~Wang, Y.~Wang, and J.~Cai.
\newblock A conservative {F}ourier pseudo-spectral method for the nonlinear
  {S}chr{\"o}dinger equation.
\newblock {\em Journal of Computational Physics}, 328:354--370, 2017.

\bibitem{GK2007}
T.~Grava and C.~Klein.
\newblock Numerical solution of the small dispersion limit of {K}orteweg-de
  {V}ries and {W}hitham equations.
\newblock {\em Communications on Pure and Applied Mathematics: A Journal Issued
  by the Courant Institute of Mathematical Sciences}, 60(11):1623--1664, 2007.

\bibitem{GK2012}
T.~Grava and C.~Klein.
\newblock A numerical study of the small dispersion limit of the {K}orteweg-de
  {V}ries equation and asymptotic solutions.
\newblock {\em Physica D: Nonlinear Phenomena}, 241(23-24):2246--2264, 2012.

\bibitem{JW92}
R.~James and J.~Weideman.
\newblock Pseudospectral methods for the {B}enjamin-{O}no equation.
\newblock {\em Advances in Computer Methods for partial differential
  equations}, 7:371--377, 1992.

\bibitem{KLM_98}
D.~J. Kaup, T.~I. Lakoba, and Y.~Matsuno.
\newblock Complete integrability of the {B}enjamin-{O}no equation by means of
  action-angle variables.
\newblock {\em Phys. Lett. A}, 238(2-3):123--133, 1998.

\bibitem{KM_98}
D.~J. Kaup and Y.~Matsuno.
\newblock The inverse scattering transform for the {B}enjamin-{O}no equation.
\newblock {\em Stud. Appl. Math.}, 101(1):73--98, 1998.

\bibitem{KK2003}
C.~Kenig and K.~Koenig.
\newblock On the local well-posedness of the {B}enjamin-{O}no and modified
  {B}enjamin-{O}no equations.
\newblock {\em Math. Res. Lett.}, 10:879--895, 2003.

\bibitem{KM2009}
C.~Kenig and Y.~Martel.
\newblock Asymptotic stability of solitons for the {B}enjamin-{O}no equation.
\newblock {\em Revista Matematica Iberoamericana}, 25:909--970, 2009.

\bibitem{KPV1994}
C.~E. Kenig, G.~Ponce, and L.~Vega.
\newblock On the generalized {B}enjamin-{O}no equation.
\newblock {\em Trans. Amer. Math. Soc.}, 342:155--172, 1994.

\bibitem{KS2015}
C.~Klein and J.-C.Saut.
\newblock {IST} versus {PDE}, a comparative study, in {H}amiltonian {P}artial
  {D}ifferential {E}quations and {A}pplications.
\newblock {\em Fields Institute Communications}, 75:383--449, 2015.

\bibitem{KP2015}
C.~Klein and R.~Peter.
\newblock Numerical study of blow-up and dispersive shocks in solutions to
  generalized {K}orteweg-de {V}ries equations.
\newblock {\em Phys. D}, 304/305:52--78, 2015.

\bibitem{KT2003}
H.~Koch and N.~Tzvetkov.
\newblock On the local well-posedness of the {B}enjamin-{O}no equation on
  {$H^s(\mathbb R)$}.
\newblock {\em Int. Math. Res. Not.}, 26:1449--1464, 2003.

\bibitem{LY2007}
T.~Lakoba and J.~Yang.
\newblock A generalized {P}etviashvili iteration method for scalar and vector
  hamiltonian equations with arbitrary form of nonlinearity.
\newblock {\em Journal of Computational Physics}, 226(2):1668--1692, 2007.

\bibitem{LRY2021}
O.~Landoulsi, S.~Roudenko, and K.~Yang.
\newblock Interaction with an obstacle in the 2d focusing nonlinear
  {S}chr\"odinger equation.
\newblock {\em arXiv preprint arXiv:2102.02170}, 2021.

\bibitem{LP2019}
U.~Le and D.~E. Pelinovsky.
\newblock Convergence of {P}etviashvili's {M}ethod near {P}eriodic {W}aves in
  the {F}ractional {K}orteweg--de {V}ries {E}quation.
\newblock {\em SIAM Journal on Mathematical Analysis}, 51(4):2850--2883, 2019.

\bibitem{LY2006}
H.~Liu and J.~Yan.
\newblock A local discontinuous {G}alerkin method for the {K}orteweg-de {V}ries
  equation with boundary effect.
\newblock {\em Journal of Computational Physics}, 215(1):197--218, 2006.

\bibitem{LY2016}
H.~Liu and N.~Yi.
\newblock A {H}amiltonian preserving discontinuous {G}alerkin method for the
  generalized {K}orteweg--de {V}ries equation.
\newblock {\em Journal of Computational Physics}, 321:776--796, 2016.

\bibitem{LZQS2019}
Z.~Liu, H.~Zhang, X.~Qian, and S.~Song.
\newblock Mass and energy conservative high order diagonally implicit
  {R}unge-{K}utta schemes for nonlinear {S}chr\"odinger equation in one and two
  dimensions.
\newblock {\em arXiv preprint arXiv:1910.13700}, 2019.

\bibitem{MRRY2021}
A.~Millet, A.~D. Rodriguez, S.~Roudenko, and K.~Yang.
\newblock Behavior of solutions to the 1d focusing stochastic nonlinear
  {S}chr{\"o}dinger equation with spatially correlated noise.
\newblock {\em Stochastics and Partial Differential Equations: Analysis and
  Computations}, pages 1--50, 2021.

\bibitem{MPST93}
T.~Miloh, M.~Prestin, L.~Shtilman, and M.~Tulin.
\newblock A note on the numerical and {N}-soliton solutions of the
  {B}enjamin-{O}no evolution equation.
\newblock {\em Wave Motion}, 17(1):1--10, 1993.

\bibitem{MR2004b}
L.~Molinet and F.~Ribaud.
\newblock Well-posedness results for the generalized {B}enjamin-{O}no equation
  with arbitrary large initial data.
\newblock {\em Int. Math. Res. Not.}, 70:3757--3795, 2004.

\bibitem{MR2004a}
L.~Molinet and F.~Ribaud.
\newblock Well-posedness results for the generalized {B}enjamin-{O}no equation
  with small initial data.
\newblock {\em J. Math. Pures Appl.}, 83:277--311, 2004.

\bibitem{NakaAk}
A.~Nakamura.
\newblock {A Direct Method of Calculating Periodic Wave Solutions to Nonlinear
  Evolution Equations. I. Exact Two-Periodic Wave Solution}.
\newblock {\em Journal of the Physical Society of Japan}, 47(5):1701--1705,
  1979.

\bibitem{OSSS2016}
D.~Olson, S.~Shukla, G.~Simpson, and D.~Spirn.
\newblock Petviashvilli's method for the {D}irichlet problem.
\newblock {\em J. Sci. Comput.}, 66(1):296--320, 2016.

\bibitem{Ono_1975}
H.~Ono.
\newblock Algebraic solitary waves in stratified fluids.
\newblock {\em Journal of the Physical Society of Japan}, 39(4):1082--1091,
  1975.

\bibitem{PS2004}
D.~E. Pelinovsky and Y.~A. Stepanyants.
\newblock Convergence of {P}etviashvili's iteration method for numerical
  approximation of stationary solutions of nonlinear wave equations.
\newblock {\em SIAM J. Numer. Anal.}, 42(3):1110--1127, 2004.

\bibitem{PD00}
B.~Pelloni and V.~Dougalis.
\newblock Numerical solution of some nonlocal, nonlinear dispersive wave
  equations.
\newblock {\em J. Nonlinear Sci.}, 10(1):1--22, 2000.

\bibitem{RWY2021}
S.~Roudenko, Z.~Wang, and K.~Yang.
\newblock Dynamics of solutions in the generalized {B}enjamin-{O}no equation: A
  numerical study.
\newblock {\em Journal of Computational Physics}, 445:110570, 2021.

\bibitem{Sanz1988}
J.~Sanz-Serna.
\newblock Runge-{K}utta schemes for {H}amiltonian systems.
\newblock {\em BIT Numerical Mathematics}, 28(4):877--883, 1988.

\bibitem{SA1991}
J.~Sanz-Serna and L.~Abia.
\newblock Order conditions for canonical {R}unge-{K}utta schemes.
\newblock {\em SIAM Journal on Numerical Analysis}, 28(4):1081--1096, 1991.

\bibitem{S2018}
J.~Saut.
\newblock Benjamin-{O}no and {I}ntermediate {L}ong {W}ave equations:
  {M}odeling, {IST} and {PDE}.
\newblock {\em In: Miller P., Perry P., Saut JC., Sulem C. (eds) Nonlinear
  Dispersive Partial Differential Equations and Inverse Scattering. Fields
  Institute Communications}, 83. Springer, New York, NY.:95--160, 2018.

\bibitem{S1979}
J.-C. Saut.
\newblock Sur quelques g\'en\'eralisations de l’ \'equation de {K}orteweg-de
  {V}ries.
\newblock {\em J. Math. Pures Appl.}, 58:21--61, 1979.

\bibitem{STL2011}
J.~Shen, T.~Tang, and L.-L. Wang.
\newblock {\em Spectral methods}, volume~41 of {\em Springer Series in
  Computational Mathematics}.
\newblock Springer, Heidelberg, 2011.
\newblock Algorithms, analysis and applications.

\bibitem{SXY2018}
J.~Shen, J.~Xu, and J.~Yang.
\newblock The scalar auxiliary variable ({SAV}) approach for gradient flows.
\newblock {\em Journal of Computational Physics}, 353:407--416, 2018.

\bibitem{SXY2019}
J.~Shen, J.~Xu, and J.~Yang.
\newblock A new class of efficient and robust energy stable schemes for
  gradient flows.
\newblock {\em SIAM Review}, 61(3):474--506, 2019.

\bibitem{T2004}
T.~Tao.
\newblock Global well-posedness of the {B}enjamin-{O}no in {$H^1(\mathbb R)$}.
\newblock {\em J. Hyperbolic Diff. Eq.}, 1(1):27--49, 2004.

\bibitem{TM98}
V.~Thom{\'e}e and A.~V. Murthy.
\newblock A numerical method for the {B}enjamin-{O}no equation.
\newblock {\em BIT Numerical Mathematics}, 38(3):597--611, 1998.

\bibitem{Tr2000}
L.~N. Trefethen.
\newblock {\em Spectral methods in {MATLAB}}, volume~10 of {\em Software,
  Environments, and Tools}.
\newblock Society for Industrial and Applied Mathematics (SIAM), Philadelphia,
  PA, 2000.

\bibitem{V2009}
S.~Vento.
\newblock Sharp well-posedness results for the generalized {B}enjamin-{O}no
  equation with high nonlinearity.
\newblock {\em Differ. Integr. Equ.}, 22(5-6):425--446, 2009.

\bibitem{Weideman95}
J.~A.~C. Weideman.
\newblock Computing the {H}ilbert transform on the real line.
\newblock {\em Math. Comp.}, 64(210):745--762, 1995.

\bibitem{Yang2021}
K.~Yang.
\newblock Arbitrarily high-order conservative schemes for the generalized
  {K}orteweg-de {V}ries equation.
\newblock {\em arXiv preprint arXiv:2103.13608}, 2021.

\bibitem{YHL2013}
N.~Yi, Y.~Huang, and H.~Liu.
\newblock A direct discontinuous {G}alerkin method for the generalized
  {K}orteweg-de {V}ries equation: energy conservation and boundary effect.
\newblock {\em Journal of Computational Physics}, 242:351--366, 2013.

\end{thebibliography}

\end{document}